\documentclass[11pt, twoside, letter]{article}

\usepackage[left = 1.1in, right = 1.1in, top = 1.5in, bottom=1in]{geometry} 
\usepackage[utf8]{inputenc}

\usepackage{hyperref}

\usepackage[english]{babel}

\usepackage{amsmath}
\usepackage{mathtools}
\usepackage{amsfonts} 
\usepackage{mathrsfs} 
\usepackage{enumitem} 
\usepackage{amssymb}
\usepackage{bbm} 
\usepackage{amsthm} 
\usepackage{pifont} 

\usepackage{nccmath} 
\let\oldint\int 
\renewcommand{\int}{\medint\oldint}

\usepackage{fancyhdr}
\renewcommand{\headrulewidth}{1pt}
\renewcommand{\headrule}{\hbox to\headwidth{%
  \color{blue4}\leaders\hrule height \headrulewidth\hfill}}

\fancypagestyle{firststyle}
{
   \fancyhf{}
   \fancyhead[R]{\thepage}
\renewcommand{\headrulewidth}{1pt}
\renewcommand{\headrule}{\hbox to\headwidth{%
  \color{blue4}\leaders\hrule height \headrulewidth\hfill}}
}

\usepackage{newtxtext,newtxmath}\selectfont
\usepackage[defaultsans]{cantarell}
\usepackage{eucal} 
\usepackage{mathtools} 

\usepackage[T1]{fontenc} 

\usepackage{titlesec, blindtext}
\usepackage[usenames,dvipsnames]{xcolor}
\definecolor{blue2}{RGB}{188,200,248} 
\definecolor{blue3}{RGB}{162,171,215}
\definecolor{blue4}{RGB}{110,113,149}
\definecolor{gris}{RGB}{153,153,153}
\definecolor{morado}{RGB}{137,104,205}
\definecolor{azulrey}{RGB}{0,0,205}
\definecolor{cornflowerblue}{RGB}{100,149,237}
\definecolor{salmon}{RGB}{250,128,114}
\definecolor{skyblue}{RGB}{135,206,250}
\definecolor{verde}{RGB}{69,139,116}
\definecolor{rosa}{RGB}{238 ,99 ,99} 
 \definecolor{gray75}{gray}{0.75}
\newcommand{\hsp}{\hspace{10pt}}
\titleformat{\section}[hang]{\Large\bfseries\color{blue4}}{\thesection\hsp\textcolor{blue4}{\textbf{$\big\bracevert$}}\hsp}{0pt}{\Large\bfseries\color{blue4}}
\titleformat{\subsection}[hang]{\large\bfseries\color{blue4}}{\thesubsection\textcolor{blue4}{\textbf{ $||$ }}}{0pt}{\large\bfseries\color{blue4}}

\usepackage[symbol,bottom,hang,flushmargin]{footmisc}

\usepackage{listings} 

\usepackage{graphicx}
\usepackage{tikz} 
\usepackage{float}
\usepackage{subcaption} 
\usepackage{wrapfig} 
\usepackage[font=small,labelfont={color=blue4,bf}]{caption} 
\usepackage{caption,setspace} 
\graphicspath{ {C:/Users/danie/Desktop/Graficos/} }

\newtheorem{theorem}{\color{blue4}Theorem \color{blue4}}[section]
\newtheorem{corollary}[theorem]{\color{blue4}Corollary \color{blue4}}
\newtheorem{lemma}[theorem]{\color{blue4}Lemma \color{blue4}}
\newtheorem{remark}[theorem]{\color{blue4}Remark \color{blue4}}

\newtheorem{proposition}[theorem]{\color{blue4}Proposition \color{blue4}}
\numberwithin{equation}{section} 
\newtagform{bbrackets}[\textbf]{\color{blue4}\textbf{(}}{\textbf{)}}
\usetagform{bbrackets}

\let\OLDthebibliography\thebibliography
\renewcommand\thebibliography[1]{
  \OLDthebibliography{#1}
  \setlength{\parskip}{2pt}
  \setlength{\itemsep}{1.4pt plus 0ex}
}

\newcommand{\mxi}{m_1}
\newcommand{\g}{\gamma}

\newcommand{\OO}{\mathcal{O}}
\newcommand{\LL}{\mathcal{L}}
\newcommand{\mxiq}{m_{q}}

\newcommand{\tS}{\tilde{S}}
\newcommand{\hatS}{\hat{S}}
\newcommand{\checkS}{\check{S}}

\newcommand{\tX}{\tilde{X}}
\newcommand{\X}{\mathcal{X}}
\newcommand{\T}{\mathbb{T}}
\newcommand{\TT}{\mathcal{T}}
\newcommand{\p}{\mathbb{P}}
\newcommand{\Q}{\mathbb{Q}}
\newcommand{\C}{\mathbb{C}}
\newcommand{\R}{\mathbb{R}}

\newcommand{\I}{\mathcal{I}}
\newcommand{\N}{\mathbb{N}}

\newcommand{\W}{\mathcal{W}}
\newcommand{\U}{\mathcal{U}}
\newcommand{\M}{\mathcal{M}}
\newcommand{\E}{\mathbb{E}}
\newcommand{\e}{\text{\normalfont e}}

\newcommand{\F}{\mathcal{F}}
\newcommand{\ZZ}{\mathcal{Z}}

\newcommand{\Var}{\text{\normalfont Var}}

\newcommand{\1}{\mathbbm{1}}
\newcommand{\s}{\sigma}
\newcommand{\dd}{\text{\normalfont d}}
\newcommand{\Int}{\text{\normalfont Int}}
\newcommand{\Leb}{\text{\normalfont Leb}}

\newcommand{\indep}{\perp \!\!\! \perp}

\newcommand\myeqd{\stackrel{\mathclap{\normalfont\mbox{\tiny{(d)}}}}{=}}

\DeclareMathOperator\supp{supp}

\usepackage{environ}
\makeatletter
\newsavebox{\measure@tikzpicture}
\NewEnviron{scaletikzpicturetowidth}[1]{%
  \def\tikz@width{#1}%
  \begin{lrbox}{\measure@tikzpicture}%
  \BODY
  \end{lrbox}%
  \pgfmathparse{#1/\wd\measure@tikzpicture}%
  \BODY
}
\makeatother

\begin{document}

\thispagestyle{firststyle}

\begin{center}

{\LARGE{\textcolor{white}{.}

Random walks with echoed steps I}}

{\large{Daniela Portillo del Valle\footnote[1]{\href{mailto:daniela.pdv@math.uzh.ch}{\nolinkurl{daniela.pdv@math.uzh.ch}}. Research supported by the Swiss National Science Foundation (SNSF), Project 212115.}}}

\vspace{-2mm}

\textit{Institute of Mathematics, University of Zürich}\footnote[2]{Winterthurerstrasse 190, 8057 Zürich, Switzerland.}\\

\end{center}


\renewcommand{\abstractname}{\begin{large}\textcolor{blue4}{Abstract}\end{large}}
\begin{abstract} A random walk with echoed steps (RWES) is a process $\{\tS_n\}_{n\geq1}=\{\tX_1+\cdots+\tX_n\}_{n\geq1}$ that inserts memory and echo into an ordinary random walk (ORW) with i.i.d. steps, $X_1+\cdots+X_n$. The RWES is defined recursively as follows. Let $\tS_1=X_1$. With probability $1-p$, the $n$-th increment of the RWES follows that of the ORW, $\tX_n=X_n$. Otherwise, $\tX_n$ is set as a random \textit{echo} of a uniform sample of the past steps $\tX_1,\dots,\tX_{n-1}$ determined by a random factor $\xi_n$. Namely, $\tX_n=\xi_n\tX_{\U[n]}$ with probability $p$, where $\U[n]\sim$Uniform$\{1,\dots,n-1\}$. The RWES is a broad generalization of the elephant random walk and of the positively/negatively/unbalanced step-reinforced random walks. We determine strong convergences of $\tS$ when the echo law $\xi$ is non-negative. The rates of convergence are determined by the product $p\E\xi$ and exhibit a phase transition with critical value at $p\E\xi=1$. Highlight that in its super-critical regime, the RWES has super-linear scaling exponents --observed for the first time in this type of random walks with memory--. We provide Laws of Large Numbers, conditions for the convergence of $\tS$ around its mean towards random series and provide some distributional properties of the limits. Our approach relies on the interpretation of the model in terms of continuous time branching random walks, random recursive trees, Pólya urns, and associated martingales.\\


\end{abstract}


\section{Introduction}\label{section:intro}

Schütz and Trimper introduced the \textit{elephant random walk} (ERW) in \cite{Schutz} to investigate how the insertion of memory in the dynamics of the ordinary random walk (ORW) may induce a super-diffusive behavior in the long-term. An elephant starts its walk at 0 and makes its first step at random in $\{-1,1\}$. At every subsequent step $n$ of its path, he remembers one of the first $n-1$ steps uniformly at random and decides to make it in the same direction with probability $p\in[0,1]$ or in the opposite direction, with probability $1-p$. The resulting process, although still inhomogeneous Markovian, exhibits an anomalous diffusive state in its convergence, as concluded by Baur and Bertoin in \cite{ERWPolyaTypeUrns} and predicted by others (c.f. \cite{daSilva}, \cite{Schutz}). Furthermore, there is a phase transition in the asymptotic behavior of the walk in terms of the memory parameter $p$, disclosing a diffusive and super-diffusive behavior towards Gaussian and non-Gaussian limits \cite{ERWPolyaTypeUrns}, respectively. We highlight that the elephant random walk has been extensively researched in the past two decades, see for example \cite{Erich}, \cite{Bercu}, \cite{Bercu2}, \cite{MultiERW}, \cite{Coletti2}, \cite{Coletti1}, \cite{daSilva}, \cite{Lucille}, \cite{Lucille2} and \cite{Kursten}.

As pointed out originally by Kürsten \cite{Kursten}, the ERW has an equivalent formulation in terms of the introduction of innovation in the system, which in turn naturally motivates generalizations of the model such as the step-reinforced random walks:\vspace{-5mm}
\begin{itemize}[leftmargin=0.5cm, label=\textcolor{blue4}{$\blacktriangleright$}]
	\item The notion of (positively) \textit{step-reinforced random walk} (SRRW), under the framework of memory insertion, was introduced by Bertoin in \cite{Noise}, \cite{ScalingExponents} and \cite{NRBM}. In this setting, the walk $\hat{S}$ starts at a point $\hat{X}_1=X_1$ with arbitrary distribution. For $n\geq2$, the $n$-th increment of the walk $\hat{X}_n$ is constructed recursively: with probability $1-p$ the walk \textit{innovates} into a new independent direction $\hat{X}_n=X_n$ with the same distribution as $X_1$. Whereas, with probability $p$, a step from the past $\hat{X}_1,\dots,\hat{X}_{n-1}$ is sampled uniformly at random and \textit{reinforced} (repeated) by the walker, i.e., $\hat{X}_n=\hat{X}_{\U[n]}$ where $\U[n]\sim$Uniform$\{1,\dots,n-1\}$.\vspace{-2mm}
	\item The introduction of negative reinforcement was formulated in \cite{Count}, also by Bertoin, under the name of \textit{counterbalanced-step random walks} (CSRW). The construction of such process $\check{S}$ follows that of the SRRW, but when the walker makes a step from the past, he incorporates it into his path going towards the opposite direction. Namely, letting the $n$-th increment of the walk to be $\check{X}_n=-\hat{X}_n$ in such case.\vspace{-2mm}
	\item A more recent approach relies on the \textit{unbalanced step-reinforced random walk} (USRRW) introduced in \cite{Aguech2}, see also \cite{HuDong}. In this setting when the walker does not innovate, he remembers a step from the past and makes it in the same direction with probability $q\in[0,1]$, following the increment of the SRRW $\hat{X}_n$, or in the opposite direction with probability $1-q$, i.e., following that of the CSRW $\check{X}_n$.
\end{itemize}
\vspace{-5mm} After the introduction and study of the SRRW and CSRW by Bertoin, these models were further investigated in \cite{Marco}, \cite{MarcoAle}, \cite{Hu}, \cite{HuZhang} and \cite{Qin1}. In these collection of works one can find an extensive study of the longtime behavior of the processes in the sense of laws of large numbers, scaling exponents, weak convergence of fluctuations and joint invariance principles.

Note from the former model descriptions, that the ERW, SRRW, CSRW and USRRW follow the same memory algorithm: sampling uniformly a step from the past with probability $p\in[0,1]$. In the same fashion, the \textit{random walk with echoed steps} (RWES) $\{\tS_n\}_{n\geq1}=\{\tX_1+\cdots+\tX_n\}_{n\geq1}$ inserts memory into an ORW $\{X_1+\cdots+X_n\}_{n\geq1}$, but also random resonance in its steps. The RWES starts at $\tS_1=X_1$. The $n$-th increment of the walk is set as the original $n$-th step of the ORW with probability $p\in[0,1]$, i.e., $\tX_n=X_n$. Whilst with probability $1-p$, $\tX_n$ is set as a random \textit{echo} of one of the steps made in the past $\tX_1,\dots,\tX_{n-1}$ chosen uniformly at random by letting $\tX_n=\xi_n\tX_{\U[n]}$, where $\xi_n$ is random, $\U[n]\sim$Uniform$\{1,\dots,n-1\}$ and $\xi_n$ and $\U[n]$ are independent. The RWES thus broadly generalizes the fore-mentioned models. When $X\equiv1$, $p=1$ and $\xi$ has Rademacher distribution, the process coincides with the ERW. Whereas when $X$ is arbitrary and $p\in(0,1]$, the RWES coincides with: \textit{i}) the SRRW if $\xi\equiv1$, \textit{ii}) the CSRW if $\xi\equiv-1$, and \textit{iii}) the USRRW if $\xi\sim$Rademacher.

In this paper we direct our efforts towards the analysis of the strong convergences of the RWES when the echo law $\xi$ is non-negative  --the general case and weak convergences will be covered in an upcoming work--. Namely, we state Laws of Large Numbers and the convergence of $\tS$ around its mean towards random series. The methods applied to this end rely on the interpretation of the model in terms of continuous time branching random walks, random recursive trees, Pólya urns and Bernoulli-bond percolation. In general, a remarkable feature about this kind of random walks with memory relies on their connection with the aforementioned models. This is why this area has aroused the interest of several authors in recent times and there are several variations of random walks with memory in the literature. To mention a few: SRRW on $\R^d$ \cite{Marco}, \cite{Qin1}, the monkey random walk \cite{MRW}, the shark random swim \cite{SharkRS}, random walks with reinforced memory of preferential attachment type \cite{Erich}, amnesic SRRW \cite{MarcoLucille}, the two-elephant walking model \cite{Aguech1}, regularly varying SRRW \cite{Majumdar1}, or the SRRW under memory lapses \cite{memorylapses}.
 

Let us assume that $X,\xi\in L^1(\p)$ and $\xi\geq0$ almost surely. Write
\[
	m_1=\E\xi.
\]
The rate of convergence of the RWES turns out to be determined by the product $p\E\xi=pm_1$ and discloses a phase transition with critical value at $pm_1=1$. This is hinted by the estimates of the expectation of the position of the walk after $n$ steps. By definition of the increments of the random walk, 
\[
	\E\tS_{n+1}=(1-p)\E X+\left(1+\frac{pm_1}{n}\right)\E\tS_n.
\]
As we shall see in Section 2, this recursive relation can be readily solved to observe that, when the memory parameter $p\in(0,1)$, 
\begin{equation}\label{eq:comportamientoasintoticoesperanzas}
	\E\tS_{n}\sim\begin{dcases}
		n^{pm_1}\times\frac{\E X}{\Gamma(1+p m_1)}\left[1+\frac{1-p}{pm_1-1}\right] &\text{ if }pm_1>1,\\
		n\log n\times(1-p) \E X &\text{ if }pm_1=1,\\
		n\times\frac{(1-p) \E X}{1-pm_1} &\text{ if }pm_1<1,
	\end{dcases}\quad\text{as }n\to\infty.\vspace{-2mm}
\end{equation}

The degenerate cases $\p(X=0)=1$, or $\p(\xi=0)=1$, or $p=0$ are omitted from this work. If $p=0$, the walk only innovates, meaning that $\tS$ is nothing more than an ORW. Similarly, if $\xi\equiv0$, $\tS$ is a random walk with independent increments with law $(1-\varepsilon)X$. Lastly, the process is identically zero if $X\equiv0$. Therefore, from this point on it will be assumed throughout the rest of this work that 
\begin{equation}\label{eq:hpt}
	p\in(0,1],\quad \p(X=0)<1,\quad \p(\xi=0)<1,\quad\xi\geq0\text{ a.s.}\quad\text{ and }\quad X,\xi\in L^\alpha(\p)\text{ for some }\alpha>1.\tag{H}
\end{equation}
For $\gamma\in\R$, let us write
\[
	m_\gamma=\E\xi^\gamma\quad\text{and}\quad\Lambda=\{\g>1:\E|X|^\g<+\infty\text{ and }m_\g<\g m_1\}.
\]

We proceed to assert our main results, starting with the super-critical regime $pm_1>1$. It is worth pointing out that, in this case, the walk $\tS$ exhibits a super-linear growth, which is a trait observed for the first time in the scheme of this type of random walks with memory. 


\begin{theorem}\label{thm:convergencesupercritical} If $pm_1>1$, then 
\[
	\lim_{n\to\infty}\frac{\tS_{n}}{n^{pm_1}}=M_\infty^{(p,\xi)}\quad\text{a.s.},
\]
where $M_\infty^{(p,\xi)}\in L^1(\p)$. Furthermore:
\vspace{-.5cm}\begin{enumerate}[leftmargin=0.5cm]
	\item[\textbf{{\color{blue4}{i)}}}] If $\E[\xi\log\xi]< m_1$, then the convergence above holds in $L^\gamma(\p)$ for every $\gamma\in\Lambda$. Moreover, the limit $M_\infty^{(p,\xi)}$ is non-degenerate and
\[
	\E M_\infty^{(p,\xi)}=\frac{\E X}{\Gamma(1+pm_1)}\left[1+\frac{1-p}{pm_1-1}\right].
\]
	\vspace{-5mm}\item[\textbf{{\color{blue4}{ii)}}}] If $\E[\xi\log\xi]\geq m_1$, then $M_\infty^{(p,\xi)}=0$ almost surely.
\end{enumerate}
\end{theorem}

By non-degenerate we mean that the limit is not constant, i.e., actually random. In sharp contrast to the ERW, SRRW, CSRW and USRRW, the convergence of the RWES is not only determined by the regime of $pm_1$, but also by a feature about the distribution of the echo variables $\xi$, namely, whether $\E[\xi\log\xi]< m_1$. This condition shall look familiar to the reader with knowledge about branching random walks and is indeed linked to the uniform integrability of martingales associated to the interpretation of the RWES in terms of the former. Understanding this condition is the key for determining the cases of convergence of the walk $\tS$ towards a non-degenerate limit. This is done in Section \ref{section:casep=1} analysing the case of pure echoing $p=1$, in which the RWES does not innovate and randomly echoes its first step forever. We will find that Theorem \ref{thm:convergencesupercritical} also holds whenever $p=1$ and $\p(\xi=1)<1$ (otherwise, $\tS_n=n X_1$), regardless of the regime of $m_1$. In addition, we characterize the condition
\[
	 \E[\xi\log\xi]<m_1\quad\text{if and only if}\quad\Lambda\text{ is non-empty}.
\]

As can be seen from the estimates (\ref{eq:comportamientoasintoticoesperanzas}), if $pm_1\leq1$, the random walk is expected to grow at a higher rate than $n^{pm_1}$. Namely, $n$ in the subcritical regime $pm_1<1$, and $n\log(n)$ in the point of criticalilty $pm=1$. In this cases, a Law of Large Numbers for the RWES can be formulated as follows. Remark that this result generalizes Theorem 1.1 in \cite{MarcoAle}, given for the case $\xi\equiv1$.


\begin{theorem}\label{thm:LLN} If $pm_1\leq1$, there exists $\gamma>1$ such that \vspace{-5mm}
\begin{enumerate}[leftmargin=0.5cm]
	\item[\textbf{{\color{blue4}{i)}}}] If $pm_1=1$ and $\E[\xi\log\xi]<m_1$, then 
	\[
		\lim_{n\to\infty}\frac{\tS_n}{n\log(n)}=(1-p)\E X,\quad\text{a.s. and in }L^\gamma(\p).
	\]
	\vspace{-5mm}\item[\textbf{{\color{blue4}{ii)}}}]  If $pm_1<1$, then
	\[
		\lim_{n\to\infty}\frac{\tS_n}{n}=\frac{(1-p)\E X}{1-pm_1},\quad\text{a.s. and in }L^\gamma(\p).
	\]
\end{enumerate}
\end{theorem}


The asymptotic behavior in (\ref{eq:comportamientoasintoticoesperanzas}) for the subcritical and critical regimes $pm_1\leq1$ suggests also that a proper re-scaling around its mean can make the RWES converge, in a similar fashion as the supercritical regime $pm_1>1$. It is thus possible to refine Theorem \ref{thm:LLN} for the subcritical and critical cases, to obtain a convergence towards a random limit, as asserted bellow. Let us address that this result generalizes Theorem 1 in \cite{NRBM}, in which $\xi\equiv1$.

\begin{theorem}\label{thm:convergencesubcritical} 
If $pm_1\leq1$ and $\E[\xi\log\xi]<m_1$, then the following assertions hold:\vspace{-5mm}
\begin{enumerate}[leftmargin=0.5cm]
	\item[\textbf{{\color{blue4}{i)}}}] If $pm_1=1$, 
	\[
		\lim_{n\to\infty}\frac{1}{n}\left[\tS_n-(1-p)\E X\text{ }n\log(n)\right]=M_\infty^{(p,\xi)}\quad\text{a.s. and in }L^\g(\p)\text{ for }\gamma\in\Lambda.
	\]
	\vspace{-5mm}\item[\textbf{{\color{blue4}{ii)}}}] If $pm_1\in(1/2,1)$ and $\Lambda\cap(1/pm_1,+\infty)\neq\emptyset$, then
	\[
		\lim_{n\to\infty}\frac{1}{n^{pm_1}}\left[\tS_n-\frac{(1-p)\E X}{1-pm_1}\text{ }n\right]=M_\infty^{(p,\xi)}\quad\text{a.s. and in }L^\g(\p)\text{ for }\gamma\in\Lambda.
	\]
\end{enumerate}
\vspace{-0.5cm}In both scenarios, the limit $M_\infty^{(p,\xi)}\in L^\gamma(\p)$ is non-degenerate and centered.
\end{theorem}

Remark that the convergence for the subcritical regime is readily seen to fail if $pm_1<1/2$. By independence, the amount of steps of the RWES in which there is only innovation is of order $n$ and, thus, the convergence of the sum of such steps $\sum_{i=1}^n(1-\varepsilon_i)X_i$ around its mean is of order $\sqrt{n}$. Thus, the rate of convergence of the re-scaling of $\tS$ is at least $1/2$. Nevertheless, the study of the regime $pm_1\in(0,1/2]$ is still interesting and will be covered in an upcoming work. It is also worth pointing out that sufficient and necessary conditions for the hypothesis in Theorem \ref{thm:convergencesubcritical}ii) are provided in Subsection \ref{subsection:subcriticalcase} below.


The methods applied to prove the results rely on three approaches, organized in this work as follows. In Section \ref{section:preliminaries} we provide a description of the model. Section \ref{section:casep=1} is devoted to prove Theorem \ref{thm:convergencesupercritical} in the pure echoing case $p=1$, and relies on the interpretation of the model through a continuous time branching random walk. This part is inspired by Baur and Bertoin \cite{Fragmentation}, and Businger \cite{SharkRS}; also we adapt arguments of Uchiyama \cite{Uchiyama} and Lyons, Pemantle and Peres \cite{Peres} to this setting. Subsequently, in Section \ref{section:trees} we find a component representation of the RWES and extend our findings in Section \ref{section:casep=1} to these components when $p\in(0,1)$. This is done linking the memory algorithm of the RWES with Bernoulli bond percolation on random recursive trees, and allows us to identify the limiting random variables $M_\infty^{(p,\xi)}$ as random series in a fashion simile to \cite{Count}, \cite{Noise} and \cite{SharkRS}. Section \ref{section:proofThms} is devoted to prove Theorems \ref{thm:convergencesupercritical}, \ref{thm:LLN} and \ref{thm:convergencesubcritical} by means of a remarkable martingale, an adaptation of that of Bercu \cite{Bercu} for the ERW, generalized by Bertenghi and Rosales-Ortiz \cite{MarcoAle} for SRRW and CSRW. Lastly, some distributional properties associated to the limits are given in Section \ref{section:propsLo}.

We conclude the introduction by pointing out that this model could find an application in behavioral economics. Namely, in the studies about the impact of the memory of the consumer on decision making. See, for example, \cite{Justinas2} and \cite{Justinas1}. In this kind of theory the premise is that the reaction to the past available information, such as former tendencies in stock markets, affects the long-time behavior of the consumer. This, under the hypothesis that there is over-reaction or under-reaction to new information depending on previous experience, thus attenuating of amplifying the original reaction, as modeled by the RWES.


\section{Preliminaries}\label{section:preliminaries}

The random walk with echoed steps (RWES) is the resulting process after incorporation of memory and resonance of steps in an ordinary random walk (ORW) $\{S_n\}_{n\geq1}$ with increments $\{X_n\}_{n\geq1}$. Let the walk start in the same random position $\tS_1=\tX_1=X_1$ as the ORW and fix a memory parameter $p\in(0,1]$. For the upcoming times $n\geq2$, we let the RWES make the same step as the ORW with probability $p$, letting $\tX_n=X_n$. Whereas, with probability $1-p$ we let the $n$-th step of the RWES be an echo of a step already made in the past by setting $\tX_n=\xi_n\tX_{\U[n]}$, where $\xi_n$ is a random quantity and $\U[n]$ a uniform sample of $\{1,\dots,n-1\}$. Specifically, the increments of the random walk with echoed steps $\tS=\{\tS_n\}_{n\geq1}$ with a fixed memory parameter $p\in(0,1]$ can be expressed as
\[
	\tX_1=X_1\quad\text{and}\quad\tX_{n}=(1-\varepsilon_{n})X_{n}+\varepsilon_{n}\xi_{n}\tX_{\U[n]},\text{ for }n\geq2,
\]
where the following are mutually independent random elements:\vspace{-.6cm}
\begin{itemize}[leftmargin=0.5cm, label=\textcolor{blue4}{$\blacktriangleright$}]
	\item A sequence of independent random variables $\{X_n\}_{n\geq1}$ with distribution $X$, that shall be referred to as the \textit{spins} of the random walk,\vspace{-2mm}
	\item $\{\varepsilon_n\}_{n\geq2}$ independent with distribution Bernoulli$(p)$,\vspace{-2mm}
	\item A collection of independent random variables $\{\U[n]\}_{n\geq2}$ such that $\U[n]\sim$Uniform$\{1,\dots,n-1\}$, and\vspace{-2mm}
	\item The \textit{echo variables}, $\{\xi_n\}_{n\geq2}$, independent and identically distributed according to the \textit{echo law} $\xi$.
\end{itemize}
\vspace{-5mm}The position of the \textit{random walk with echoed steps} $\tS=\{\tS_n\}_{n\geq1}$ at time $n\geq1$ is given by 
\[
	\tS_n=\tX_1+\cdots+\tX_n.
\]
The name \textit{random walks with echoed steps} was chosen eluding the fact that the introduction of the random weight $\xi$ into the reinforcement appears in the model with a multiplicative effect of the previous steps in the form of an \textit{echo}. Each time there is reinforcement in the path of the walk, the process can remember, not only the initial step $\tX_1$ or the innovation steps in which $\tX_n=X_n$, but also the formerly reinforced steps already resonated by at least one random quantity $\xi$. This way, when the process reinforces, a reflection, amplification or attenuation of the remembered step is added to the path of the walker. 

A transparent way to visualize the model is through labeled uniform random recursive trees (URRT), a path already explored for other step-reinforced random walks in \cite{Fragmentation}, \cite{MarcoAle}, \cite{Noise}, \cite{ScalingExponents}, \cite{NRBM}, \cite{Count}, \cite{SharkRS}, to mention a few. An URRT is an increasing tree that is constructed recursively as follows. The tree starts with a root 1 and at time $n=2$ a leaf is attached to the root with label 2. For the upcoming steps $n\geq3$, one of the $n$ vertices of the tree is chosen uniformly and independently at random and the $(n+1)$-th leaf is attached to this randomly chosen vertex, say, to $\U[n+1]\sim$Uniform$\{1,\dots,n\}$.

\begin{figure}
	\centering
	\begin{subfigure}{0.4\linewidth}
		\centering
		\begin{tikzpicture} [
		    every node/.style={circle,solid, draw=black,thick, minimum size = 0.5cm, line width=0.5mm},
		    emph/.style={edge from parent/.style={dashed, red, draw, line width=0.5mm}},
		    norm/.style={edge from parent/.style={solid, black, draw, line width=0.5mm}}
		    ]
		  \node [draw, black, circle] {1}
		    [black, line width=0.5mm]
		    child { 
		        node [draw, black, circle] {2} 
		        child [solid] { 
		            node [draw, black, circle] {4} 
		            child { node [draw, black, circle] {5} }
		        }
		        child {
		            node [draw, black, circle] {6} 
		            child { node [draw, black, circle] {7} }
		            child { node [draw, black, circle] {8} 
		                child { node [draw, black, circle] {9} }
		            }
		        }
		    }
		    child { node [draw, black, circle] {3} }
		    child { node [draw, black, circle] {10} 
		        child { node [draw, black, circle] {11} }
		    }
		    ;
		    
		\end{tikzpicture}
		\caption{$p=1$.}
	\end{subfigure}
	\hspace{1cm}
	\begin{subfigure}{0.4\linewidth}
		\centering
		
		\begin{tikzpicture} [
		    every node/.style={circle,solid, draw=black,thick, minimum size = 0.5cm, line width=0.5mm},
		    emph/.style={edge from parent/.style={dashed,red,draw,line width=0.5mm}},
		    norm/.style={edge from parent/.style={solid,black,line width=0.5mm,draw}}
		    ]
		  \node [draw, rosa, circle] {1}
		    [rosa, line width=0.5mm]
		    child [white, dashed, line width=0.25mm]{
		        node [draw, blue4, circle] {2} 
		        child [solid, blue4, line width=0.5mm] { 
		            node [draw, blue4, circle] {4} 
		            child { node [draw, blue4, circle] {5} }
		        }
		        child [solid,blue4, line width=0.5mm] {
		            node [draw, blue4, circle] {6} 
		            child { node [draw, blue4, circle] {7} }
		            child [white, dashed, line width=0.25mm] { node [draw, verde, circle] {8} 
		                child [solid, verde, line width=0.5mm] { node [draw, verde, circle] {9} }
		            }
		        }
		    }
		    child { node [draw, rosa, circle] {3} }
		    child { node [draw, rosa, circle] {10} 
		        child { node [draw, rosa, circle] {11} }
		    }
		    ;
		\end{tikzpicture}
		\caption{$p\in(0,1)$.}
	\end{subfigure}
	\put(-210,135){\color{gris}{\LARGE$\rightsquigarrow$}}
	\put(-330,182.5){\small$\xi_2$}
	\put(-325,140){\small$\xi_6$}
	\put(-305,97.5){\small$\xi_8$}
	\put(-292.5,50){\small$\xi_9$}
	\put(-76.125,50){\color{verde}{\small$\xi_9$}}
	\caption{Samples of the memory tree associated to the random walk with echoed steps $\tS$ for $p=1$ and $p\in(0,1)$.}
	\label{fig:arbolito}
	\vspace{-.5cm}
\end{figure}
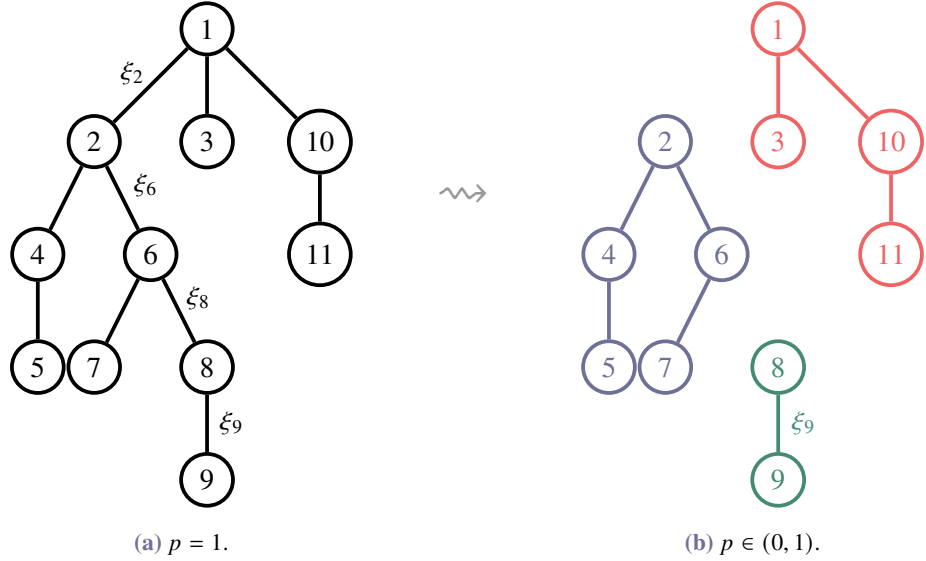

To have a better idea of the behavior of the RWES, $\tS$, consider the case without innovation $p=1$. Then the dynamic of $\tS$ consists in adding random echoes of its first step $\tS_1=X_1$. We associate the initial step of the walk $\tX_1=X_1$ to the root of the tree 1. At $n=2$, the only step the walker can remember is the initial step, thus $\U[2]=1$, $\tX_2=\xi_2X_1$ and the vertex 2 is attached to the root 1. For the step $n=3$, the walker can remember the steps 1 or 2, which happens according to the random variable $\U[3]\sim$Uniform$\{1,2\}$. If $\U[3]=1$, the third step of the walk is $\tX_3=\xi_3 X_1$ and we attach the vertex 3 to the root 1. Whereas if $\U[3]=2$, $\tX_3=\xi_1\xi_2 X_1$ and vertex 3 is attached to 2. This construction follows recursively, so the \textit{memory tree} of $\tS$ is a URRT. Figure \ref{fig:arbolito}(a) displays a realization of such random construction. 

Every time the RWES remembers a step from the past, it is echoed through a random quantity with the law of $\xi$, all independent. We thus associate to the edges of the memory tree one random weight corresponding to the echo variable, $\xi_n$ is placed on the edge connecting vertex $n$ to its parent $\U[n]$. The value of $\tX_k$ for $k\in\{1,\dots,n\}$ is then given by locating the node $k$ in the tree and multiplying all the weights $\xi_i$ on the edges in the path up to the root 1 of the tree. For example, in Figure \ref{fig:arbolito} we have $\tX_9=\xi_2\xi_6\xi_8\xi_9X_1$. Thus the position of the walk after $n$ steps coincides with the sum of the multiplicative weights of all the nodes going up to the root of the tree, multiplied by $X_1$.

Now let $p\in(0,1)$. After the construction of the recursive tree as for $p=1$, remove each one of the edges independently with probability $1-p$, i.e., perform Bernoulli bond percolation on the memory URRT. This is equivalent to the innovation of the walk $\tS$ according to the variables $\{1-\varepsilon_n\}_{n\geq2}$. As before, the value of the increment $\tX_k$ is found multiplying the random weights $\xi$ up the way to the root of the sub-tree to which it belongs after the partition, say the vertex $r$, and the variable $X_r$. In Figure \ref{fig:arbolito}(b), we have $\varepsilon_8=0$ and now the 9-th increment of the walk $\tS$ with $p\in(0,1)$ is given simply by $\tX_9=\xi_9X_8$. 

Let us start by giving explicit computations of the expected value of the random walk $\tS$. Recall the assumptions (\ref{eq:hpt}) and write
\[
	m_q=\E\xi^q,\quad\text{for }q\in\R.
\]
Note that, due to the independence structure of $(\varepsilon,\xi, X,\U)$, for every non-negative (or bounded) measurable function $f$,
\[
	\E f(\tX_{n+1})=(1-p)\E f(X)+\frac{p}{n}\sum_{k=1}^n\E f(\xi_{n+1}\tX_k).
\]
This identity gives the following recurrence relation for the sums of the moments of the increments $\tX_k$ letting $f(x) = |x|^q$:
\begin{equation}\label{eq:recursionmomentos}
	\sum_{k=1}^{n+1}\E|\tX_k|^q=(1-p) \E|X|^q+\frac{n+pm_q}{n}\sum_{k=1}^n\E|\tX_k|^q,
\end{equation}
where the absolute value can be dropped if $q=1$. Applying (\ref{eq:recursionmomentos}) inductively yields the computation of the expected values in the lemma below.


\begin{lemma}\label{lemma:recursionmedias} Let $q\in\R$. If $X,\xi\in L^q(\p)$, then
		\[
			\sum_{k=1}^{n}\E|\tX_k|^q=\begin{dcases}
				\frac{\E|X|^q}{1-p\mxiq}\left[(1-p)n+\frac{p(1-\mxiq)}{\Gamma(1+p\mxiq)}\frac{\Gamma(n+pm_q)}{(n-1)!}\right]&\quad\text{if }pm_q\neq1,\\
				n\E|X|^q\left[p+(1-p)\sum_{j=1}^{n}\frac{1}{j}\right]&\quad\text{if }pm_q=1,
			\end{dcases}
		\]
		where the absolute value can be dropped for $q=1$.
\end{lemma}

From Lemma \ref{lemma:recursionmedias} one gets the estimates for $\E\tS_n$ as $n\to\infty$ stated in (\ref{eq:comportamientoasintoticoesperanzas}) given that as $n\to\infty$, $\sum_{j=1}^nj^{-1}\sim\log(n)$ and that, by Stirling's approximation,
\begin{equation}\label{eq:stirling}
	\frac{\Gamma(n+a)}{\Gamma(n+b)}\sim n^{a-b}.
\end{equation}
These estimates suggest on first hand the rates of convergence of $\tS$ in function of the product $pm_1$. More generally, from Lemma \ref{lemma:recursionmedias} one also gets that $\sum_{i=1}^n\E|\tX_i|^q\sim C n^{1\vee(pm_q)}$ if $pm_q\neq1$ and, thus $pm_q$ has a phase transition with critical point at $pm_q=1$, where $\sum_{i=1}^n\E|\tX_i|^q\sim C n\log(n)$. 

We conclude the preliminaries by remarking a property of the gamma function that will be frequently applied and which proof readily follows by induction:
\begin{equation}\label{eq:sumagammas}
	\sum_{i=m}^n\frac{\Gamma(i+a)}{\Gamma(i+1+b)}=\frac{1}{b-a}\left[\frac{\Gamma(m+a)}{\Gamma(m+b)}-\frac{\Gamma(n+1+a)}{\Gamma(n+1+b)}\right] \quad\text{for }1\leq m\leq n\text{ and }a\neq b.
\end{equation}


\section{The case of pure echoing: $p=1$}\label{section:casep=1}

Our way to prove the convergences asserted in Theorems \ref{thm:convergencesupercritical} and \ref{thm:convergencesubcritical}, starts with the analysis of the random walk that echoes its initial step $X_1\equiv1$ all along its path with probability 1. That is, throughout this section it will assumed that the memory parameter and the initial position of the walk are
\[
	p=1\quad\text{and}\quad X_1\equiv1.
\]
This way, the walk never innovates into a new direction and its memory tree is never cut. In this case, we can write directly
\[
	\tX_1=1,\quad\tX_{n+1}=\xi_{n+1}\tX_{\U[n+1]}\quad\text{and}\quad\tS_{n}=\sum_{k=1}^n\tX_k.
\]
Note also that if $p=1$, the assumption $X_1\equiv1$ results in no loss of generality. Recall the notation
\[
	m_\gamma=\E\xi^\gamma\quad\text{and}\quad\Lambda=\{\g>1:\E|X|^\g<+\infty\text{ and }m_\g<\g m_1\}.
\]
The main statement of the section is the following. 

\begin{proposition} \label{prop:convergencecasepeq1} If $p=1$ and $X\equiv1$, then
\[
	\lim_{n\to\infty}\frac{\tS_{n}}{n^{m_1}}=\LL^{(\xi)}\quad\text{almost surely},
\]
where $\LL^{(\xi)}\in L^1(\p)$. Furthermore:\vspace{-0.5cm}
\begin{enumerate}[leftmargin=0.5cm]
	\item[\textbf{{\color{blue4}{i)}}}] If $\E[\xi\log\xi]< m_1$, then the convergence above holds in $L^\gamma(\p)$ for every $\gamma\in\Lambda$. Moreover, the limit $\LL^{(\xi)}\in L^\gamma(\p)$ is non-degenerate if $\p(\xi=1)<1$ and
\[
	\E\LL^{(\xi)}=\frac{1}{\Gamma(1+ m_1)}.
\]
	\vspace{-5mm}\item[\textbf{{\color{blue4}{ii)}}}] If $\E[\xi\log\xi]\geq m_1$, then $\LL^{(\xi)}=0$ almost surely.
\end{enumerate}
\end{proposition}

The proofs of Proposition \ref{prop:convergencecasepeq1}i) and \ref{prop:convergencecasepeq1}ii) require different methods. The first one, fairly elementary, relies on finding an associated martingale and showing that it is bounded in $L^\g(\p)$ under the conditions of i). This can be found in Subsection \ref{subsection:||T||acotLa}. Moreover, we will find that, in fact,
\[
	\E[\xi\log\xi]<m_1\quad\text{if and only if}\quad\Lambda_0:=\{\g>1:m_\g<\g m_1\}\neq\emptyset.
\]
For the proof of ii), in Subsection \ref{subsection:BRW} we define a branching random walk (BRW) in continuous time and an associated exponential additive martingale, which convergence is tightly connected to that of the pure echoing RWES. Subsequently, in Subsection \ref{subsection:proofWoo=0} we state a ``many to one formula'' for the BRW and provide a Biggins' convergence-like theorem to give proof of Proposition \ref{prop:convergencecasepeq1}ii).


\subsection{Proof of Proposition \ref{prop:convergencecasepeq1}i)}\label{subsection:||T||acotLa}

To give proof of the first part of Proposition \ref{prop:convergencecasepeq1}i), we will make use of a convergent martingale associated to the RWES, inspired in that of \cite{Bercu} and \cite{MarcoAle} in their studies about the ERW, SRRW and CSRW, respectively. Then we shall derive the non-degeneracy of the limit in virtue of the following proposition.


\begin{proposition}\label{prop:||T(n)||ismtg} Let $p=1$ and $X\equiv1$. The process 
\[
	M_n=\frac{(n-1)!}{\Gamma(n+ m_1)} \tS_n,\quad n\geq1
\]
is a martingale with respect to the filtration $\F$, where $\F_n=\s(\{\xi_k,\U[k]\}_2^n)$ for $n\geq2$ and $\F_1$ is trivial. In addition, if $\xi\in L^\gamma(\p)$ for some $\gamma>1$ such that $m_\gamma<\gamma m_1$, then the martingale is bounded in $L^\gamma(\p)$.
\end{proposition}


Taking Proposition \ref{prop:||T(n)||ismtg} for granted, $M$ converges almost surely towards an integrable limit that we shall denote by $\LL^{(\xi)}$. Then, by Stirling's approximation,
\[
	\lim_{n\to\infty}M_n=\lim_{n\to\infty}\frac{\tS_n}{n^{m_1}}=\LL^{(\xi)} \quad\text{almost surely}.
\]
Furthemore, if $\xi\in L^\g(\p)$ for some $\g>1$ such that $m_\gamma<\gamma m_1$, then $M$ is bounded in $L^\gamma(\p)$, and hence $M$ converges to $\LL^{(\xi)}$ in $L^\gamma(\p)$ as well. This ensures that $\E\LL^{(\xi)}=\E M_1=\Gamma(1+ m_1)^{-1}$. Furthermore, $\LL^{(\xi)}$ is non-degenerate whenever $\p(\xi=1)<1$, since otherwise $\tS_n/n\equiv1$. Therefore, Proposition \ref{prop:convergencecasepeq1}i) follows from Proposition \ref{prop:||T(n)||ismtg} addressing that $\E[\xi\log\xi]< m_1$ is equivalent to $m_\gamma<\gamma m_1$ for some $\gamma>1$. This is stated properly in the following lemma, which proof is elementary and postponed to Appendix \ref{section:appendixUchiyama}.


\begin{lemma}\label{lemma:derivadamtheta} Assume that $\p(\xi=1)<1$. Let $\I=\Int\{a\geq0:m_a<+\infty\}$ and $\phi_\theta(r)=m_{r\theta}/r$ for fixed $\theta\in \I$. The mapping $\phi_\theta$ is twice differentiable and strictly convex on $\I$, and attains its minimum at $r_\theta\in\I$ if and only if $\phi_\theta(r_\theta)=\theta\E[\xi^{\theta r_\theta}\log\xi]$. Furthermore,
\[
	\theta\E[\xi^\theta\log\xi]<m_\theta\quad\text{if and only if}\quad\text{there exists }\gamma>1\text{ such that } m_{\gamma\theta}<\gamma m_\theta.
\]
\end{lemma}


We move over to the proof of Proposition \ref{prop:||T(n)||ismtg}. This involves the application of the following inequalities to show the $L^\gamma(\p)$-boundedness of the martingale $M$, which proof can be found in Appendix \ref{section:appendixUchiyama}.

\begin{lemma}\label{lemma:inequalities} Fix $\gamma>2$ and $\delta>0$. There exist constants $N_{\gamma,\delta},N_{\gamma,\delta}'>0$ such that for every $z,w\in\R$,
\[
	|z+w|^\gamma\leq\begin{dcases}
		|w|^\gamma+\gamma|w|^{\gamma-2}wz+N_{\gamma,\delta}|w|^{\gamma-1}|z|&\quad\text{ if }|z|\leq\delta|w|,\\
		|w|^\gamma+\gamma|w|^{\gamma-2}wz+N_{\gamma,\delta}'|z|^\gamma&\quad\text{ if }|z|>\delta|w|, 
	\end{dcases}
\]
where
\[
	\lim_{\delta\downarrow0^+}N_{\gamma,\delta}=0\quad\text{and}\quad\lim_{\delta\downarrow\infty}N_{\gamma,\delta}'=1.
\]
If $\gamma\in(1,2]$ the second inequality holds for every $w,z\in\R$ regardless of $|z|>\delta|w|$ and $N_{\gamma,\delta}'$ does not depend on $\gamma$.
\end{lemma}


The last ingredient required to prove Proposition \ref{prop:||T(n)||ismtg} is a distributional equation for $\tS$ under the current assumptions, $p=1$ and $X_1\equiv1$. For this, we exploit the interpretation of uniform random recursive trees in terms of Pólya urns (c.f. \cite{Chauvin} and \cite{Lucille} for similar arguments). 

\begin{lemma}\label{lemma:fixedpointeq||T||} If $p=1$ and $X\equiv1$, then
\[
	\{\tS_{n+1}\}_{n\geq1}\myeqd \{\hatS_{R_n}+\xi\checkS_{B_n}\}_{n\geq1},
\]
where $ \hatS,\checkS$ are copies of $\tS$, $\{(R_n,B_n)\}_{n\geq1}$ is a Pólya-Eggenberger urn with initial condition $(R_1,B_1)=(1,1)$ and deterministic replacement matrix $I_2$, and all the elements above are mutually independent.
\end{lemma}
\vspace{-1cm}\begin{proof}[\bf\color{blue4}\textit{Proof.}] Recall that the memory tree of the RWES is a uniform random recursive tree. Consider $\{(R_n,B_n)\}_{n\geq1}$ a classical two color Pólya urn with initial condition $(R_1,B_1)=(1,1)$ and replacement matrix $I_2$, where $R_n$ (resp. $B_n$) denotes the number of red balls (resp. blue balls) in the urn after the $n$-th step, i.e., after $n-1$ extractions. At the second step of evolution of $\tS$, the vertex 2 is attached to the root 1 of its memory tree and from then on the following vertices are successively attached to a vertex chosen uniformly at random. Thus, after the $(n+1)$-th step of the walk, the amount of vertices in the tree that do not have the vertex 2 as an ancestor has the same distribution as the amount of red balls in the urn after $n-1$ extractions, $R_{n}$. While the amount of vertices that have the vertex 2 as an ancestor has the same distribution as the amount of blue balls after $n-1$ extractions, $B_n=n+1-R_n$. This means that right after its $(n+1)$-th step, the amount of times that the RWES has randomly echoed the step $\tX_1=\xi_1$ is distributed as $n+1-R_n$, whereas that of the rest of the steps is distributed as $R_n$. The distributional equation then follows recalling that $\{\xi_n\}_{n\geq2}$ and $\{\U[n]\}_{n\geq2}$ are independent.
\end{proof}


\begin{proof}[\bf\color{blue4}\textit{Proof of Proposition~\ref{prop:||T(n)||ismtg}}] First, since $\tX_{n+1}=\xi_{n+1}\tX_{\U[n+1]}$ and $\{(\xi_{n},\U[n])\}_{n\geq2}$ are independent,
\[
	\E\left[\tS_{n+1}\big|\{\xi_i,\U[i]\}_{i=2}^n\right]=\tS_n+\sum_{i=1}^n\E\left[\xi_{n+1}\tX_i\1_{U[n+1]=i}\Big|\{\xi_i,\U[i]\}_{i=2}^n\right]=\left(\frac{n+m_1}{n}\right)\tS_n;
\]
where we have used that $\U[n+1]$ is a uniform random variable on $\{1,\dots,n\}$. It follows that $M$ is a martingale by multiplying both sides by $n!/\Gamma(n+1+ m_1)$.

Suppose there exists $\gamma>1$ such that $m_\gamma<\gamma m_1$. We shall prove that $M$ is bounded in $L^\gamma(\p)$ by showing inductively over $\lfloor\gamma\rfloor\geq1$ that, under such condition,
\begin{equation}\label{eq:aux||T||boundedLgamma}
	\text{there exists }C_\gamma>0\text{ such that for every }n\geq1,\quad\E \tS_n^\gamma\leq C_\gamma \frac{\Gamma(n+\gamma m_1)}{(n-1)!}.
\end{equation}
This suffices in virtue of the asymptotic (\ref{eq:sumagammas}). First, assuming that $\gamma\in(1,2]$ and letting $w=\tS_n$ and $z=\tX_{n+1}$ in Lemma \ref{lemma:inequalities}ii),
\[
	\E\left[\tS_{n+1}^\gamma\big|\{\xi_i,\U[i]\}_{i=2}^n\right]\leq\left(1+\frac{\gamma m_1}{n}\right)\tS_n^\gamma+N\frac{m_\gamma}{n}\sum_{i=1}^n\tX_i^\gamma.
\]
By Lemma \ref{lemma:recursionmedias} and (\ref{eq:stirling}), taking expectations in the last display yields
\[
	\E\tS_{n+1}^\gamma\leq\left(\frac{n+\gamma m_1}{n}\right)\E\tS_n^\gamma+N \frac{\Gamma(n+m_\gamma)}{n!\Gamma(m_\gamma)}.
\]
Applying the latter inequality inductively and expressing the resulting products in terms of gamma functions, gives the following in virtue of (\ref{eq:stirling}),
\[\begin{split}
	\E\tS_{n+1}^\gamma\leq\frac{\Gamma(n+1+\gamma m_1)}{n!}\times\frac{1}{\Gamma(1+\gamma m_1)}\left[1+N\frac{m_\gamma}{\gamma m_1-m_\gamma}\right].
\end{split}\]

On the other hand, if $\gamma>2$, according to Lemma \ref{lemma:inequalities} we can choose $\delta=\delta(\gamma)$ sufficiently large so that $N_{\gamma,\delta(\gamma)}'<\gamma m_1/m_\gamma$. Note that Lemma \ref{lemma:derivadamtheta} guarantees that if $m_\gamma<\gamma m_1$, then also $m_{\gamma'}<\gamma'm_1$ for every $\gamma'\in(1,\gamma)$ due to the convexity of $\phi_1$. So we can define $C_\gamma>0$ recursively as follows
\[
	C_\gamma=\begin{dcases}
		\frac{2}{\Gamma(1+\gamma m_1)}\left\{1\vee \left[\frac{(\gamma+N_{\gamma,\delta})m_1}{\gamma m_1-m_\gamma N_{\gamma,\delta}'}C_{\gamma-1}\Gamma(1+(\gamma-1) m_1)\right]\right\}&\quad\text{if }\gamma>2,\\
		\frac{1}{\Gamma(1+\gamma m_1)}\left[1+N\frac{m_\gamma}{\gamma m_1-m_\gamma}\right]&\quad\text{if }\gamma\in(1,2].
	\end{dcases}
\]
It was previously showed that (\ref{eq:aux||T||boundedLgamma}) holds for this defintion of $C_\gamma$ in the case $\gamma\in(1,2]$. We now proceed to argue that (\ref{eq:aux||T||boundedLgamma}) is true for $\gamma>2$, assuming inductively that it holds for every $\gamma'\in(1,\lfloor\gamma\rfloor]$. The case $n=1$ is trivial since $X_1\equiv1$ and $C_\gamma\Gamma(1+\gamma m_1)\geq 2$ by definition. For $n\geq2$, writing the addends of the distributional equation of Lemma \ref{lemma:fixedpointeq||T||}, as $w=\hatS_{U}$ and $z=\xi\checkS_{n-U}$ in Lemma \ref{lemma:inequalities}ii),
\[
	\E\tS_{n}^\gamma=\E\left[(\hatS_{U}+\xi\checkS_{n-U})^\gamma\right]\leq\frac{1+N'_{\gamma,\delta}m_\gamma}{n-1}\sum_{i=1}^{n-1}\E\tS_i^\gamma+\frac{(\gamma+N_{\gamma,\delta})m_1}{n-1}\sum_{i=1}^{n-1}\E\tS_{i}^{\gamma-1}\E\tS_{n-i},
\]
where we have used that $\hatS_{U}\myeqd\checkS_{n-U}$, since $U$ is uniformly distributed over $\{1,\dots,n-1\}$. Applying Lemma \ref{lemma:recursionmedias} to $\E\tS_{n-i}$ and (\ref{eq:aux||T||boundedLgamma}) to $\E\tS_i^{\gamma-1}$ gives
\[
	\sum_{i=1}^{n-1}\E\tS_i^{\gamma-1}\E\tS_{n-i-1}\leq C_{\gamma-1}\frac{\Gamma(n+\gamma m_1)}{(n-2)!}\frac{\Gamma(1+(\gamma-1)m_1)}{\Gamma(2+\gamma m_1)},
\]
where we have completed the Beta function $\beta(i+(\gamma-1)m_1,n-i+m_1)$ inside the resulting sum and iterated with the integral sign to obtain this expression. Thus, by applying the induction hypothesis (\ref{eq:aux||T||boundedLgamma}) to $\E\tS_i^\gamma$ for $i\in\{1,\dots,n-1\}$ and solving the remaining sum with (\ref{eq:sumagammas}),
\[\begin{split}
	\E\tS_{n}^\gamma
	&\leq\frac{\Gamma(n+\gamma m_1)}{(n-1)!}\left[C_\gamma\frac{1+N'_{\gamma,\delta}m_\gamma}{1+\gamma m_1}+(\gamma+N_{\gamma,\delta})m_1 C_{\gamma-1}\frac{\Gamma(1+(\gamma-1)m_1)}{\Gamma(2+\gamma m_1)}\right]\leq\frac{\Gamma(n+\gamma m_1)}{(n-1)!}C_\gamma,
\end{split}\]
where the last inequality holds by definition of $C_\gamma$. Then, (\ref{eq:aux||T||boundedLgamma}) holds for $\gamma>2$, completing the proof.
\end{proof}


\begin{remark} \normalfont Note that it is not true in general that if $\xi\in L^\alpha(\p)$ and $\E[\xi\log\xi]<m_1$, then the stated convergences hold in $L^\alpha(\p)$. Proposition \ref{prop:||T(n)||ismtg} (and therefore Proposition \ref{prop:convergencecasepeq1}i) for $p=1$) assures that there exists $\gamma>1$, possibly $\gamma<\alpha$, such that the convergence holds in $L^\gamma(\p)$. For example, if $\xi$ is a standard exponential random variable, we have that $\xi\in L^2(\p)$ and $m_\gamma/\gamma=\Gamma(\gamma)$. Thus, the function $\phi_1(\gamma)=m_\gamma/\gamma$ reaches its minimum at a point $r^*\approx1.4616$ and $m_\gamma/\gamma<1=m_1$ for every $\gamma\in(1,2)$. Then, Proposition \ref{prop:convergencecasepeq1}i) ensures that $n^{-1}\tS_n\to\LL^{(\xi)}$ in $L^{\gamma}(\p)$ for $\gamma\in(1,2)$. Nevertheless, it does not guarantee that also in $L^2(\p)$. Regardless of this, by Lemma \ref{lemma:derivadamtheta} and Proposition \ref{prop:convergencecasepeq1}i), there always exists $\theta\in(0,1)$ such that the convergence holds both a.s. and in $L^{\g}(\p)$ for $\tS_n^{(\theta)}=\sum_{i=1}^n\tX_i^\theta$ , the RWES with echoing variables $\{\xi^\theta_n\}_{n\geq2}$.
\end{remark}


\subsection{A continuous time branching random walk}\label{subsection:BRW}

So far, the condition $\E[\xi\log\xi]<m_1$ has been proven to be sufficient for the limit $\LL^{(\xi)}$ of the re-scaled RWES to be non-degenerate when the memory parameter $p=1$ and the initial step of the walk is $X_1\equiv1$. We turn our to attention to show that such condition is, in fact, necessary and complete the proof of Proposition \ref{prop:convergencecasepeq1}. As previously mentioned, a different approach is carried out to this end. Namely, exploiting the connection between the RWES and a class of continuous time branching random walks, which are introduced in this section and that will be found quite useful further on for the analysis of the RWES in the general case $p\in(0,1]$.

Consider under $\p$ a branching process made of static particles located in $\R$ starting with an initial common ancestor located at position 0 and that evolves according to the following dynamics. Particles reproduce independently, and every alive particle gives birth to new born particles (one at the time) at independent exponential times with parameter $1$. The position of the new born is the one of its father, displaced by an independent random variable with law $\log\xi\in[-\infty,\infty)$. The particles never die and, in particular, the number of particles alive at time $t\geq0$ is a Yule process $(Y_t)_{t\geq0}$ with rate $1$. That is, a pure-birth process which starts at 1, and with birth rate from state $n\in\N$ given by $n$. With a slight abuse of notation we still write $Y_t$ for the set of particles alive at time $t\geq0$. 

Denote the point measure whose atoms consist of the particles alive in the system at time $t\geq0$ as
\[
	\ZZ_t(\dd x)=\sum_{x\in Y_t}\delta_x(\dd x),
\]
and write $(\F_t)_{t\geq0}$ for its canonical filtration. It is plain from the construction that $\ZZ=(\ZZ_t(\dd x))_{t\geq0}$ is a continuous time Markov chain taking values in the set of finite measures on $\R$, which is referred to as a continuous time branching random walk. This process traces back to \cite{Biggins}, \cite{Ikeda}, \cite{Watanabe}, \cite{Uchiyama} and was recently further studied in \cite{JeanHairuo}, to mention a few. We refer to this references therein for details. The Markov property of $\ZZ$ reads as follows. Conditionally on $\F_t$, the process $(\ZZ_{t+s}(\dd x))_{s\geq0}$ is distributed as the sum of $Y_t$-independent copies of $\ZZ$, each of which has atoms respectively shifted by each $x\in Y_t$. This property is commonly referred to as the \textit{branching property} and also holds in the strong Markov sense, a detailed formal treaty about it can be found in \cite{Ikeda}.

The link between $\ZZ$ and $\tS$ lies in random recursive trees. Given that the particles $\{0,\dots,n\}$ are alive, the $(n+1)$-th particle is the son of $i\in\{0,\dots,n\}$ if and only if the exponential clock of particle $i$ ticks first than the others in the interval $(\tau_{n},+\infty)$, where $\tau_n=\inf\{t\geq0:Y_t\geq n\}$ denotes the time of birth of the $n$-th particle. Due to the strong Markov property, conditioned on $\tau_{n}$, the first clock to tick in $(\tau_{n},+\infty)$ starting with the particle 0 in the system, has the same distribution as the first clock to tick in $\R^+$ when the system starts with the $n+1$ particles $\{0,\dots,n\}$. Recalling that the clocks are independent exponential variables, the probability that the particle $i\in\{0,\dots,n\}$ gives birth to the particle $n+1$ is therefore uniform over $\{0,\dots,n\}$. Thus, the genealogical tree of the particles in $\ZZ_{\tau_n}$ is a uniform random recursive tree, so as the memory tree of $\tS_n$. Moreover, since $p=1$ and $X_1\equiv1$, we can write
\[
	\tS_n=1+\sum_{k=2}^{n}\e^{\log\xi_{k}+\log\tX_{\U[k]}}.
\]
Thus, we deduce that for every $\theta\geq0$, there is the identity in law
\[
	\left\{\sum_{x\in Y_{\tau_{n-1}}}\e^{\theta x}\right\}_{n\geq1}\myeqd\left\{\tX_1^\theta+\cdots+\tX_{n}^\theta\right\}_{n\geq1}.
\]
In particular, when $\theta=1$ the process in the left-hand side of the last display is a version of the RWES $\tS$. Therefore, the asymptotic behaviour of $\tS$ is intimately related with that of
\[
	\Sigma_t^{(\theta)}=\sum_{x\in Y_t}\e^{\theta x}
\]
as $t\to\infty$. In this direction, the following martingale will be crucial in our study. Recall the notation
\[
	m_\theta=\E\xi^\theta.
\]


\begin{lemma}\label{lemma:Wesmtg} Fix $\theta\geq0$ such that $\xi\in L^\theta(\p)$ and let
\[
	\W^{(\theta)}_t=\e^{-m_\theta t}\Sigma_t^{(\theta)}.
\] 
The process $\W^{(\theta)}=(\W^{(\theta)}_t)_{t\geq0}$ is a non-negative martingale with mean 1. Therefore, it converges almost surely as $t\to\infty$ to an integrable limit that will be denoted as $\W^{(\theta)}_\infty$.
\end{lemma}
\vspace{-.5cm}\begin{proof}[\bf\color{blue4}\textit{Proof.}] Let us start by showing that the mean of $\Sigma_t^{(\theta)}$ is $\e^{m_\theta t}$. Conditioning on the birth-time of the first child of the original ancestor gives 
\[\begin{split}
	\E\Sigma^{(\theta)}_{t}&=\int_0^t\left\{\E\left[\int_{\R} \e^{\theta x}\ZZ_{t-s}(\dd x)\right]+\E\left[\e^{\theta\log\xi}\int_{\R} \e^{\theta x}\ZZ_{t-s}(\dd x)\right]\right\}\e^{-s}\dd s+\int_t^\infty \e^{-s}\dd s,
\end{split}\]
where $\ZZ$ and $\xi$ are independent. Therefore, with an adequate change of variables, the latter writes as
\[
	\E\Sigma^{(\theta)}_{t}=\e^{-t}\int_0^t(1+m_\theta)\E\Sigma^{(\theta)}_u\e^{u}\dd u+\e^{-t}.
\]
Differentiating this expression with respect to $t$ gives $\frac{\dd}{\dd t}\E[\Sigma^{(\theta)}_{t}]=m_\theta\E[\Sigma^{(\theta)}_{t}]$. Solving the ordinary differential equation with initial condition $\E\Sigma^{(\theta)}_0=1$ we find
\begin{equation}\label{eq:meansigma}
	\E\Sigma^{(\theta)}_{t}=\e^{m_\theta t}.
\end{equation}

To show that $\W^{(\theta)}$ is a martingale, in virtue of the branching property at time $s<t$, we can write
\[\begin{split}
	\E\left[\sum_{x\in Y_t}\e^{\theta x}\Big|\F_s\right]=\sum_{x\in Y_s}\e^{\theta x}\E\left[\sum_{y\in Y_{t-s}}\e^{\theta y}\right]=\e^{(t-s)m_\theta}\sum_{x\in Y_s}\e^{\theta x}.
\end{split}\]
Multiplying both sides by $\e^{-m_\theta t}$ completes the proof.
\end{proof}


We proceed to give a characterization of the convergence of the martingale $\W^{(\theta)}$ in terms of the condition $\E[\xi\log\xi]<m_1$.

\begin{proposition}\label{prop:convmtgW} Fix $\theta\in\Int\{a>0:m_a<+\infty\}$. Then:\vspace{-0.5cm}
\begin{enumerate}[leftmargin=0.5cm]
	\item[\textbf{{\color{blue4}{i)}}}] If $\theta\E[\xi^\theta\log\xi]<\E\xi^\theta$, the convergence $\W^{(\theta)}_t\to\W^{(\theta)}_\infty$ holds in $L^\gamma(\p)$-sense for some $\gamma>1$,
	\item[\textbf{{\color{blue4}{ii)}}}] Otherwise, if $\theta\E[\xi^\theta\log\xi]\geq\E\xi^\theta$, then $\W^{(\theta)}_\infty=0$ almost surely.	
\end{enumerate}
\end{proposition}

It is worth pointing out that Proposition \ref{prop:convmtgW}i) for $\gamma\in(1,2]$ is covered by Theorem 5 in \cite{Biggins}, formerly proven in \cite{Uchiyama} with stronger conditions (not met by the current setting of the RWES). Thus, with Proposition \ref{prop:convmtgW}i) we extend the results of Biggins and Uchiyama in the particular setting of RWES for the case in which $\gamma>2$. 

The proof of Proposition \ref{prop:convmtgW}ii) is postponed to the next subsection, since some self involved technical results will be required. Taking it for granted, we have enough tools to prove Proposition \ref{prop:convergencecasepeq1}ii).

\begin{proof}[\bf\color{blue4}\textit{Proof of Proposition ~\ref{prop:convergencecasepeq1}ii)}] By Proposition \ref{prop:convergencecasepeq1}, $n^{-m_1}\tS_{n}$ converges almost surely to $\LL^{(\xi)}$ as $n\to\infty$. Recall that $\tau_n$ denotes the time of birth of the $n$-th particle and that $\{\Sigma_{\tau_n}^{(1)}\}_{n\geq0}$ is a version of $\tS$. Then,
\begin{equation}\label{eq:limSncomoWyE}
	\LL^{(\xi)}=\lim_{n\to\infty}\frac{\tS_n}{(n+1)^{m_1}}\myeqd\lim_{n\to\infty}\frac{\W_{\tau_n}^{(1)}}{(\e^{-\tau_n}Y_{\tau_n})^{m_1}}=\frac{\W_\infty^{(1)}}{\mathcal{E}_\infty^{m_1}},
\end{equation}
where we have used the well known fact that $\e^{-t}Y_t$ converges almost surely to a standard exponential random variable $\mathcal{E}_\infty$. Thus, $\LL^{(\xi)}\myeqd\W^{(1)}_\infty \mathcal{E}_\infty^{-m_1}$ and from Proposition \ref{prop:convmtgW}ii) follows that $\LL^{(\xi)}=0$ almost surely if $\E[\xi\log\xi]\geq m_1$.
\end{proof}


\subsection{Proof of Proposition \ref{prop:convmtgW}}\label{subsection:proofWoo=0}

Note that Proposition \ref{prop:convmtgW}i) follows from Proposition \ref{prop:||T(n)||ismtg}. Indeed, if $\theta\E[\xi^\theta\log\xi]<m_\theta$, Lemma \ref{lemma:derivadamtheta} guarantees the existence of $\gamma>1$ such that $m_{\gamma\theta}<\gamma m_\theta$. Then, Proposition \ref{prop:||T(n)||ismtg} yields that the RWES with echo variables distributed as $\xi^\theta$ is convergent almost surely and in $L^\gamma(\p)$ towards $\LL^{(\xi^\theta)}$. The non-degeneracy of $\W_\infty^{(\theta)}$ follows from noting that $\W_\infty^{(\theta)}$ with position displacements $\log\xi$ has exactly the same construction as $\W_\infty^{(1)}$ with displacement $\log(\xi^\theta)$.

We therefore direct out efforts in the remainder of this section to prove Proposition \ref{prop:convmtgW}ii). Namely, we give a Biggins' convergence-like theorem and a ``many to one formula'' for the branching random walk $\ZZ$, adapting the arguments in \cite{Peres} to the continuous time setting. Note that for every $\theta\in\Int\{a>0:m_a<+\infty\}$ there is a unique probability measure $\Q^{(\theta)}$ on $\F_\infty=\s(\cup_{t\ge0}\F_t)$ such that
\begin{equation}\label{eq:restriccionQ}
	\Q^{(\theta)}(A)=\E[\W_t^{(\theta)}\1_A],\quad t\geq0,\text{ }A\in\F_t.
\end{equation}
Its existence follows by Kolmogorov's extension theorem and the actual existence of the BRW $\ZZ$ (c.f. \cite{Ikeda}). We direct our efforts towards proving the proposition below, from which Proposition \ref{prop:convmtgW}ii) follows.

\begin{proposition}\label{prop:existenceQandconvW} For every $\theta\in\I$, the following holds:
\vspace{-0.5cm}\begin{enumerate}[leftmargin=0.5cm]
	\item[\textbf{{\color{blue4}{i)}}}] $\theta\E[\xi^\theta\log\xi]<\E\xi^\theta\quad \Leftrightarrow \quad\Q^{(\theta)}\ll\p\quad \Leftrightarrow\quad \E\W^{(\theta)}_\infty=1$,
	\item[\textbf{{\color{blue4}{ii)}}}] $\theta\E[\xi^\theta\log\xi]\geq\E\xi^\theta\quad \Leftrightarrow \quad\Q^{(\theta)}\perp\p\quad \Leftrightarrow\quad \W^{(\theta)}_\infty=0$ $\p$-a.s.
\end{enumerate}
\end{proposition}


The tool provided to show Proposition \ref{prop:existenceQandconvW} is the following ``many-to-one formula'' for the branching random walk $\ZZ$. Let us introduce some notation. For every fixed $t \geq 0$ and $x \in Y_t $ we write $(x(s):s \in [0,t])$ for the ancestral path of $x$. That is, $x(s)$ is the position of its youngest alive ancestor at time $s$, with the convention that $x$ is its own youngest alive ancestor provided it was already born. In particular, $x(t) = x$. Note that $(x(s):s \in [0,t])$ is therefore the rcll path that starts at the origin, with the original ancestor at 0, and jumps according to the accumulated positions of the ancestors of $x$ at their times of birth. Our version of the many-to-one formula writes as follows.

\begin{proposition}\label{prop:manytoone} Fix $t\geq0$ and write $D([0,t], \R)$ for the space of rcll functions from $[0,t]$ to $\R$ endowed with the Skorokhod topology. Let $\X_t^\bullet = (\X^\bullet_t(s):s \in [0,t])$ be the process which law in $D([0,t], \R)$ under $\p^\bullet$ is defined by the relation
\begin{equation}\label{eq:distXbullet}
	\E^\bullet f(\X_t^\bullet) = \e^{-t}\mathbb{E} \left[\sum_{x \in Y_t}f(x(s):s \in [0,t]) \right],
\end{equation} 
for every measurable functional $f:D([0,t], \R)\to\R_+$. Then, under $\p^\bullet$ the process $\X^\bullet_t$ is a compound Poisson process restricted to $[0,t]$ with intensity 1 and jump distribution $\log\xi$.
\end{proposition}\vspace{-.5cm}
\begin{proof}[\bf\color{blue4}\textit{Proof.}] 
Since $\E Y_t = \e^{t}$, the distribution of $\X_t^\bullet$ under $\p^\bullet$ is a probability measure. To show that $\X_t^\bullet$ has independent and stationary increments, remark that by the branching property, for every $h \geq 0$ the law of $\X_t^\bullet$ and the law of $(\X_{t+h}^\bullet(s): s \in [0,t])$ under $\mathbb{P}^\bullet$ coincide. In addition, for every $0 \leq r \leq s \leq t$ and every $Z$ non-negative $\mathcal{F}_{r}$-measurable random variable, 
    \[
         \mathbb{E}^\bullet\left[ Z f(\X_t^\bullet(s) - \X_t^\bullet(r)) \right]=\mathbb{E}^\bullet Z\mathbb{E}^\bullet \left[ f(\X_{t-r}^\bullet(s-r)) \right]=\mathbb{E}^\bullet Z\mathbb{E}^\bullet \left[ f(\X_{t}^\bullet(s-r))\right],
    \]
    in virtue of the branching property and our first observation. Therefore, under $\mathbb{P}^\bullet$,  $\X_t^\bullet$ has stationary, independent increments, and thus is a Lévy process on $[0,t]$. Moreover, $\X_t^\bullet$ is piecewise constant by construction. We thus conclude that it is a compound Poisson process under $\p^\bullet$ and its Laplace functional at time $t\geq0$ is given by
\[
	\E^\bullet \exp\{ \lambda \X_t^\bullet(t)\}=e^{-t}\E \left[\sum_{x \in Y_t}\e^{\lambda x} \right]=\e^{ t (m_\lambda - 1) }. 
\]
\end{proof}


We now establish the link between Proposition \ref{prop:manytoone} with the limiting random variable $\W_\infty^{(\theta)}$. For fixed $t\geq0$ and $\theta\in\Int\{a>0:m_a<+\infty\}$, conditionally on $\F_t$, we pick at random one of the particles $x\in Y_t$ alive at time $t$, proportionally to $\e^{\theta x}$ and write $\X_t^{(\theta)}=(\X_t^{(\theta)}(s):s\in[0,t])$ for its historical path. That is, 
\[
	\p(\X_t^{(\theta)}(t)=x|\F_t)=\frac{\e^{\theta x}}{\Sigma^{(\theta)}_t},\quad x\in Y_t.
\]
Note that by construction,
\begin{equation}\label{eq:ineqQXt}
	\e^{m_\theta t}\W_t^{(\theta)}=\sum_{x\in Y_t}\e^{\theta x}\geq\e^{\theta \X_t(t)}\quad\p|_{\F_t}\text{-a.s.}
\end{equation}
Then, in virtue of the many-to-one formula, it is possible to show that,

\begin{corollary}\label{cor:bajoQvesPPois} For fixed $t\geq0$, under $\Q^{(\theta)}$ the process $\X_t^{(\theta)}$ is a compound Poisson process restricted to $[0,t]$ with rate $m_\theta$ and jump-distribution $\Delta^{(\theta)}$ characterized by the relation 
\[
	\E f(\Delta^{(\theta)})=\frac{\E[\xi^\theta f(\log\xi)]}{m_\theta},
\]
for measurable $f:\R\to\R_+$.
\end{corollary}
\vspace{-0.5cm}\begin{proof}[\bf\color{blue4}\textit{Proof.}] Note that by construction of $\X_t^{(\theta)}$ and $\Q^{(\theta)}$, for every $g: D([0,t],\R)\to\R_+$ we have
\[
	\Q^{(\theta)} [g(\X_t^{(\theta)})]=\e^{-m_\theta t}\E\left[\sum_{x\in Y_t}g(x(s):s\in[0,t])\e^{\theta x}\right]=\frac{\E^\bullet[g(\X_t^\bullet)\exp\{\theta\X_t^\bullet(t)\}]}{\E^{\bullet}\exp\{\theta\X_t^\bullet(t)\}},
\]
where the process $\X^\bullet_t$ and the measure $\p^\bullet$ are defined in (\ref{eq:distXbullet}). Our many-to-one formula in Proposition \ref{prop:manytoone} states that $\X^\bullet_t$ under $\p^\bullet$ is a compound Poisson process restricted to $[0,t]$. Therefore, so is $\X_t^{(\theta)}$ under $\Q^{(\theta)}$, since it standard to show that the exponential tilting of a compound Poisson process is one as well. Lastly, by the last display we readily compute
\[
	\Q^{(\theta)} \exp\{\lambda\X_t^{(\theta)}(t)\}=\e^{t(m_{\theta+\lambda}-m_\theta)}=\exp\left\{m_\theta t\text{ }\E\left[\e^{\lambda\Delta^{(\theta)}}-1\right]\right\},
\]
from which the assertion follows.
\end{proof}


\begin{remark}\label{remark:manytoone}\normalfont Note that Proposition \ref{prop:manytoone} jointly with Corollary \ref{cor:bajoQvesPPois} give an expression more commonly known as many-to-one formula, namely
\[
	\E\left[\sum_{x\in {Y}_{t}}f(x(s): s \in [0,t])\right]=\E\left[\e^{-\theta P^{(\theta)}(t)+m_\theta t}f\left(P^{(\theta)}(s): s \in [0,t]\right)\right],
\]
where $\theta\in\Int\{a>0:m_a<+\infty\}$, $f:D(\mathbb{R}_+ , \R)\to\R_+$ is measurable and $P^{(\theta)}$ is a compound Poisson process with intensity $m_\theta$ and the jump distribution $\Delta^{(\theta)}$ defined in Corollary \ref{cor:bajoQvesPPois}.
\end{remark}


All concepts needed to show Proposition \ref{prop:existenceQandconvW} have been provided, so we end the section with its proof.

\begin{proof}[\bf\color{blue4}\textit{Proof of Proposition~\ref{prop:existenceQandconvW}}] First, from (\ref{eq:restriccionQ}) we have that $\Q^{(\theta)}|_{\F_t}\ll\p|_{\F_t}$. Then, by Theorem 5.3.3 in \cite{Durrett},
\begin{equation}\label{eq:durret}
	\Q^{(\theta)}(A)=\E\left[\limsup_{t\to\infty}\W_t^{(\theta)}\1_A\right]+\Q^{(\theta)}\left(A\cap\left\{\limsup_{t\to\infty}\W_t^{(\theta)}=+\infty\right\}\right),\quad\text{for every }A\in\F_\infty.
\end{equation}
Remark that $\lim_{t\to\infty}\W_t^{(\theta)}=\W_\infty^{(\theta)}$ exists $\Q^{(\theta)}$-a.s., so we can substitute the limit superior in the last display by $\W^{(\theta)}_\infty$. Indeed, from equation (\ref{eq:restriccionQ}) readily follows that $\Q^{(\theta)}[\1_A/\W_s^{(\theta)}]=\p(A)$ for every $s\in[0,t]$ and $A\in\F_s$. Then, $\E_\Q[1/\W_t^{(\theta)}|\F_s]=1/\W_s^{(\theta)}$ for every $0\leq s<t$ and $1/\W^{(\theta)}$ is a non-negative $\F_t$-martingale with respect to $\Q^{(\theta)}$. Hence, under $\Q^{(\theta)}$, $1/\W^{(\theta)}$ converges to $1/\W_\infty^{(\theta)}\in[0,\infty)$. This ensures the $\Q^{(\theta)}$-a.s. existence of $\W_\infty\in(0,+\infty]$.

On one hand, as argued at the end of last subsection, if $\theta\E[\xi^\theta\log\xi]<m_\theta$ then $\E\W_\infty^{(\theta)}=1$ by Proposition \ref{prop:||T(n)||ismtg}. Recalling that $\W_\infty^{(\theta)}<+\infty$ under $\p$, identity (\ref{eq:durret}) gives
\[
	\Q^{(\theta)}\left(\W_\infty^{(\theta)}<+\infty\right)=\E\left[\W_\infty^{(\theta)}\1_{\W_\infty<+\infty}\right]=1.
\]
Therefore (\ref{eq:durret}) writes as $\Q^{(\theta)}A=\E[\W_\infty^{(\theta)}\1_A]$ for every $A\in\F_\infty$. We conclude that $\Q^{(\theta)}\ll\p$ and $\frac{\dd\Q^{(\theta)}}{\dd\p}=\W_\infty$, in this case.

Secondly, assume that $\theta\E[\xi^\theta\log\xi]\geq m_\theta$. From expression (\ref{eq:ineqQXt}) and Corollary \ref{cor:bajoQvesPPois} follows that for every $t\geq0$ and $M>0$,
\begin{equation}\label{eq:auxQandconvW}
	\Q^{(\theta)}(\W_t^{(\theta)}\geq M)\geq\Q^{(\theta)}\left(\e^{\theta \X_t^{(\theta)}-tm_\theta}\geq M\right)=\p\left(\frac{P^{(\theta)}(t)}{t}-\frac{\log M}{\theta t}\geq\frac{m_\theta}{\theta}\right),
\end{equation}
where $P^{(\theta)}$ is a compound Poisson process with intensity $m_\theta$ and jump distribution $\Delta^{(\theta)}$. If $\theta\E[\xi^\theta\log \xi]\geq m_\theta$ then $\E P^{(\theta)}(1)=\E[\xi^\theta\log\xi]>0$, and thus, $t^{-1}P^{(\theta)}(t)\to\E[\xi^\theta\log\xi]$ $\p$-a.s in virtue of Theorem 7.2 in \cite{Kyp}. Therefore, the limit as $t\to\infty$ of right-hand side of \ref{eq:auxQandconvW} is 1 for every $M>0$. Hence, $\W_\infty^{(\theta)}=+\infty$ almost surely with respect to $\Q^{(\theta)}$ and
\[
	0=\Q^{(\theta)}(\W_\infty^{(\theta)}<+\infty)=\E\W_\infty^{(\theta)},
\]
by (\ref{eq:durret}). This in turn implies that $\W_\infty^{(\theta)}=0$, $\p$-a.s., and for every $A\in\F_\infty$,
\[
	\Q^{(\theta)}(A)=\Q^{(\theta)}(A\cap\{\W_\infty^{(\theta)}=+\infty\})\quad\text{and}\quad\p(A)=\p(A\cap\{\W_\infty^{(\theta)}<+\infty\}).
\]
Hence, $\Q^{(\theta)}\perp\p$ and the proof is complete.
\end{proof}


\section{The link with random recursive trees}\label{section:trees}

After understanding the convergence of the RWES for the pure echoing case $p=1$, we proceed to extend the strong convergences to the case $p\in(0,1)$. In that direction, we give a decomposition of $\tS$ in terms of random recursive trees that will allow us to determine the conditions of convergence and the expressions of the limits stated in Theorems \ref{thm:convergencesupercritical} and \ref{thm:convergencesubcritical}.

Recall from the preliminaries that the memory tree of the RWES is an URRT, an increasing random tree that starts with the root 1 at time 1 and to which more vertices are added successively choosing a parent uniformly and independently at random. Write $\T$ for the tree with infinite vertices and  $\T(n)$ for the tree up to time $n\geq1$. Remark that the size of the tree at time $n$ is $|\T(n)|=n$.

In the case $p=1$, the random walk echoes its initial step $\tS_1=X_1$ forever. We add the i.i.d. echo variables $\{\xi_n\}_{n\geq2}$ to the edges of the memory tree $\T$ as weights, where $\xi_n$ is placed on the edge that connects the vertex $n\geq2$ to its parent $\U[n]$. Then, we define recursively the weight $\omega(n)$ of the vertex $n$ as the that of its parent $\omega(\U[n])$ echoed by the independent variable $\xi_n$. That is, we set $\omega(n)=\xi_n\omega(\U[n])$ for $n\geq2$ and give the root weight 1, $\omega(1)=1$. Then the weight of vertex $n$ coincides with the product of the echo variables $\xi$ of all the edges connecting $n$ to the root of the tree (c.f. Figure \ref{fig:arbolito}(a)). Define the weight of the whole tree grown up to time $n$, as 
\[
	||\T(n)||=\sum_{i=1}^n\omega(i),\quad \text{where } \omega(1)=1\text{ and  }\omega(n+1)=\xi_{n+1}\omega(\U[n+1]).
\]
Then, it is plain from the construction that
\[
	\{\tS_n\}_{n\geq1}=\{X_1||\T(n)||\}_{n\geq1}.
\]

So as to introduce innovation in the model (i.e., to let the memory parameter $p\in(0,1)$), given the infinite tree $\T$ we perform Bernoulli bond percolation with parameter $p$ on its edges. That is, we destroy every edge of the tree independently at random with probability $1-p$, i.e., according to $\{1-\varepsilon_n\}_{n\geq2}$. At time $n$, this results in the original tree $\T(n)$ being broken down to at most $n$ sub-trees rooted at the vertices whose connecting edge to the parent was destroyed. Write $T_r$ for the sub-tree of $\T$ rooted at vertex $r\geq1$ after percolation, with the convention that $T_r(n)=\emptyset$ if $\varepsilon_r=1$. We proceed to weight these sub-trees by their vertices as before, according to the product of the weights $\xi$ on its edges up to the root of the respective sub-tree. Specifically, we define the weight of vertex $n\geq2$ to be 1 if its edge to its parent has been destroyed after the Bernoulli bond percolation and to be the one associated to its parent weighted by $\xi_n$, otherwise. That is,
\begin{equation}\label{eq:defwtree}
	\omega(1)=1\quad\text{and}\quad \omega(n+1)=(1-\varepsilon_{n+1})+\varepsilon_{n+1}\xi_{n+1}\omega(\U[n+1]).
\end{equation}
Thus, one can consider again the total weight of every sub-tree $T_r(n)$ at time $n\geq1$ rooted at $r\geq1$ and denote it as follows
\[
	|| T_r(n)||=\sum_{i\in T_r(n)}\omega(i),
\]
where we have written $i\in T_r(n)$ to indicate that $i$ is a vertex of the sub-tree $T_r(n)$ and with the convention $\sum_{i\in\emptyset}=0$. Note that cutting the edges of $\T$ according to $1-\varepsilon$ coincides with the event of $\tS_n$ innovating to an independent new direction. Then, if $\varepsilon_r=0$ for $2\leq r\leq n$, the total echoing in $\tS_n$ towards the direction $X_r$ coincides with $||T_r(n)||$. Thus, we can identify 
\begin{equation}\label{eq:Snarbolitos}
	\tS_{n}=\sum_{r=1}^nX_r||T_r(n)||,\quad\text{for }n\geq1.
\end{equation}
This representation is the key to determine the convergence of the RWES when the memory parameter $p\in(0,1)$. Remark that the identities in the last displays are point-wise, not distributional.

For fixed $n\geq1$, conditionally on their sizes, the sub-trees $T_1(n),\dots, T_n(n)$ that partition the memory tree of $\tS$ are independent URRT (c.f. Lemma 3.3 in \cite{Count}). Thus, in light of Proposition \ref{prop:convergencecasepeq1} , we can approach the convergence of the sub-trees $T_r$ through that of $\T$. Subsequently, we determine the convergences of $\tS$ in Theorems \ref{thm:convergencesupercritical} and \ref{thm:convergencesubcritical} via the decomposition (\ref{eq:Snarbolitos}). In this direction, let us exploit the relationship between URRT and Pólya urns to relate the distribution of $||T_r(\bullet)||$ to that of $||\T(\bullet)||$.


\begin{lemma}\label{lemma:distYnr} For every $2\leq r\leq n$, there is the equality in distribution
\begin{equation}\label{eq:eqdistTkn}
	|| T_r(n)||\myeqd(1-\varepsilon_r)|| T_1(Y(n,r))||,
\end{equation}
where $\varepsilon_r$ is Bernoulli$(p)$ distributed, $(\varepsilon_r,T_1,Y(n,r))$ are independent and
\[
	\p(Y(n,r)=i)=(r-1)\frac{(n-r)!}{(n-1)!}\frac{(n-i-1)!}{(n-r+1-i)!},
\]
for $i\in\{1,\dots,n-r+1\}$, and 0, otherwise.
\end{lemma}
\vspace{-1cm}\begin{proof}[\bf\color{blue4}\textit{Proof.}] In a similar setting as \cite{SharkRS}, for fixed $n\geq2$, after assigning a parent to the $r$-th vertex ($r\in\{2,\dots,n\}$) the tree $\T(r)$ has $r$ nodes and the rest $n-r$ nodes must be attached to it recursively. If we consider an urn with one red ball and $r-1$ blue balls (representing the $r$-th vertex of the tree and the remaining $r-1$ nodes that make up the tree, respectively), making $n-r$ extractions on the classical setting of Pólya urns coincides with the experiment of assigning parents to the remaining $n-r$ vertices in the tree. Then, the number of red ball draws in the urn after $n-r$ extractions coincides with the amount of vertices attached to the sub-tree rooted at $r$ at time $n$. This random quantity, which we denote by $Y'(n,r)$, is known to be Beta-Binomial$(n-r,1,r-1)$ distributed (c.f. Theorem 3.1 \cite{Mahmoud}), i.e.,
\[
	\p(Y'(n,r)=i)={n-r \choose i}\frac{\beta(i+1,(n-r)-i+(r-1))}{\beta(1,r-1)}=(r-1)\frac{(n-r)!(n-i-2)!}{(n-1)!(n-r-i)!},
\]
for $i\in\{0,\dots,n-r\}$ and 0, otherwise. Furthermore, if $Y'(n,r)$ is a Beta-Binomial random variable independent of $(T,\varepsilon,\xi)$, we have the following equality in distribution in virtue of the description above and the independence of the echoing variables:
\[
	|| T_r(n)||\myeqd(1-\varepsilon_r)|| T_1(1+Y'(n,r))||,
\]
where $\varepsilon_r\sim$Bernoulli$(p)$ is independent of $T_1$ and $Y'(n,r)$. The claim follows letting $Y(n,r)=1+Y'(n,r)$.
\end{proof}

As a first application of Lemma \ref{lemma:distYnr}, we compute the expected value of the weight of the sub-trees $T_r(n)$ that will be needed in our analysis.


\begin{lemma}\label{lemma:momentosT0n} For every $\alpha\in\R$ and $1\leq r\leq n$, we have
\begin{equation}\label{eq:esperanzaTrn}
	\E\left[\sum_{i\in T_r(n)}\omega(i)^\alpha\right]=(1-p\1_{r>1})\frac{(r-1)!}{(n-1)!}\frac{\Gamma(n+pm_\alpha)}{\Gamma(r+pm_\alpha)}\quad\text{for every }1\leq r \leq n.
\end{equation}
\end{lemma}
\vspace{-1cm}\begin{proof}[\bf\color{blue4}\textit{Proof.}] The equation holds for $r=1$ and $p=1$ by Lemma \ref{lemma:recursionmedias}. Now let $p\in(0,1)$ and note that 
\[
	\sum_{i\in T_1(n)}\omega(i)^\alpha\myeqd||\T^{(\varepsilon\xi^\alpha)}(n)||, 
\]
where $\T^{(\varepsilon\xi^\alpha)}$ is the memory tree associated to the RWES with echo variables $\{\varepsilon_n\xi_n^\alpha\}_{n\geq2}$ and memory parameter $1$. Then (\ref{eq:esperanzaTrn}) holds for $r=1$ by our first observation. Lastly, fix $2\leq r\leq n$. From (\ref{eq:eqdistTkn}) follows that
\[\begin{split}
	\E\left[\sum_{i\in T_r(n)}\omega(i)^\alpha\right]=(1-p)\frac{r-1}{(n-1)!}\frac{\Gamma(n+pm_\alpha)}{\Gamma(1+pm_\alpha)}\sum_{i=0}^{n-r} \binom{n-r}{i}\beta(i+1+pm_\alpha,n-i-1).
\end{split}\]
Identity (\ref{eq:esperanzaTrn}) is readily recovered by exchanging the sum and integral (of the beta function) signs in the last display.
\end{proof}


We are ready to assert the convergences of the weights in percolation components, which is in fact a corollary of Proposition \ref{prop:convergencecasepeq1}.

\begin{corollary} \label{cor:convTreesRootedatk} For every $r\geq1$,
\[
	\lim_{n\to\infty}\frac{||T_r(n)||}{n^{pm_1}}=\LL^{(p,\xi)}_r\quad\text{almost surely},
\]
where $\{\LL_r^{(p,\xi)}\}_{r\geq1}\subset L^1(\p)$ are marginally distributed as follows
\begin{equation}\label{eq:eqdistLr}
	\LL_1^{(p,\xi)}\myeqd\LL^{(\varepsilon\xi)}\quad\text{and}\quad\LL_r^{(p,\xi)}\myeqd(1-\varepsilon_r)\LL^{(\varepsilon\xi)}\beta_r^{pm_1}\quad\text{ for every }r\geq2,
\end{equation}
where $\LL^{(\varepsilon\xi)}$ is defined in Proposition \ref{prop:convergencecasepeq1}, $\varepsilon_r\sim$Bernoulli$(p)$, $\beta_r\sim$Beta$(1,r-1)$ and $(\varepsilon_r,\LL^{(\varepsilon\xi)},\beta_r)$ are mutually independent. \vspace{-0.5cm}
\begin{enumerate}[leftmargin=0.5cm]
	\item[\textbf{{\color{blue4}{i)}}}] If $\E[\xi\log\xi]< \E\xi$, then the convergences above hold in $L^\gamma(\p)$ for every $\g\in\Lambda_0=\{\g>1:m_\g<\g m_1\}$. Moreover, $\{\LL_r^{(p,\xi)}\}_{r\geq1}$ are non-degenerate and
\[
	\E\LL_r^{(p,\xi)}=(1-p\1_{r>1})\frac{(r-1)!}{\Gamma(r+pm_1)},\quad\text{for every }r\geq1.
\]
	\item[\textbf{{\color{blue4}{ii)}}}] If $\E[\xi\log\xi]\geq \E\xi$, then $\LL_r^{(p,\xi)}=0$ almost surely for every $r\geq1$.
\end{enumerate}
\end{corollary}
\vspace{-0.75cm}\begin{proof}[\bf\color{blue4}\textit{Proof.}] The case $p=1$ corresponds to Proposition \ref{prop:convergencecasepeq1} and has been already proved. Let $p\in(0,1)$ and write 
\[
	\F_1=\s(X_1),\quad\F_n=\sigma\left(X_1,\{\varepsilon_i,\xi_i,\U[i],X_i\}_{i=2}^n\right)\quad\text{and}\quad\F=\{\F_n\}_{n\geq1}.
\]
Analogously to Proposition \ref{prop:||T(n)||ismtg}, one can show that
\begin{equation}\label{eq:Tr(n)esmtg}
	\left\{\frac{(n-1)!}{\Gamma(n+pm_1)}||T_r(n)|| \right\}_{n\geq r}\text{ is a martingale with respect to }\{\F_n\}_{n\geq r},
\end{equation}
for every $r\geq1$. Then, for every $r\geq1$, we have that $n^{-pm_1}||T_r(n)||$ converges almost surely to an integrable random variable $\LL_r^{(p,\xi)}\in L^1(\p)$.

If $\T'$ is the memory tree associated to a RWES with echo variables $\varepsilon\xi$ and memory parameter $1$, then $T_1\myeqd\T'$. So the statement follows directly for $r=1$ from Proposition \ref{prop:convergencecasepeq1}. Now fix $r\geq2$. According to Theorem 3.2 in \cite{Mahmoud}, the random variable $Y(n,r)$ in (\ref{eq:eqdistTkn}) is such that $n^{-1}Y(n,r)\to\beta_r$ almost surely, where $\beta_r$ has distribution Beta$(1,r-1)$. Then, the distributional equality (\ref{eq:eqdistTkn}) gives
\[
	\LL_r^{(p,\xi)}\myeqd\lim_{n\to\infty}(1-\varepsilon_r)\frac{||T_1(Y(n,r))||}{Y(n,r)^{pm_1}}\left(\frac{Y(n,r)}{n}\right)^{pm_1}=(1-\varepsilon_r)\LL^{(\varepsilon\xi)}\beta_r^{pm_1},
\]
where $\varepsilon_r$, $\LL^{(\varepsilon\xi)}$ and $\beta_r$ are independent. If $\E[\xi\log\xi]\geq m_1$, then $\E[\varepsilon\xi\log(\varepsilon\xi)]\geq \E[\varepsilon\xi]$ and Proposition \ref{prop:convergencecasepeq1}ii) guarantees that $\LL^{(\varepsilon\xi)}=0$ and, thus, $\LL_r^{(p,\xi)}=0$ almost surely. Otherwise, if $\E[\xi\log\xi]<m_1$, since $T_1\myeqd\T'$, by Proposition \ref{prop:convergencecasepeq1}i), for every $\g\in\Lambda_0$,
\[
	\sup_{n\geq1}\left\{\frac{(n-1)!^\gamma}{\Gamma(n+pm_1)^\gamma}\E||T_1(n)||^\gamma \right\}<+\infty.
\]
Then, by (\ref{eq:eqdistTkn})
\[
	\E||T_r(n)||^\gamma=(1-p)\E||T_1(Y(n,r))||^\gamma\leq(1-p)\E||T_1(n)||^\gamma,
\]
and thus $n^{-pm_1}||T_r(n)||\to\LL_r^{(p,\xi)}$ in $L^\gamma(\p)$ as well.
\end{proof}


\section{Proof of the main results}\label{section:proofThms}

In the latter section we found a representation of $\tS$, when the memory parameter $p\in(0,1)$, in terms of its memory sub-trees
\[
	\tS_{n}=\sum_{r=1}^nX_r||T_r(n)||,\quad n\geq1
\]
and established the convergence of the components $n^{-pm_1}||T_r(n)||$ separately. This is enough to assert Theorem \ref{thm:convergencesupercritical} for the supercritical regime $pm_1>1$; while the subcritical and critical regimes $pm_1\leq1$ in Theorem \ref{thm:LLN} and \ref{thm:convergencesubcritical} require an extra element to wrap up together the convergence of $\tS$. Namely, a re-scaling of the RWES, that turns out to be a martingale. We inspect this martingale in Subsection \ref{subsection:convM} and proceed to prove the main theorems \ref{thm:convergencesupercritical}, \ref{thm:convergencesubcritical} and \ref{thm:LLN}, respectively, in the subsequent subsections.

\begin{remark}\label{remark:Lambdanonempty}\normalfont Before we proceed, let us highlight that, under the hypotheses (\ref{eq:hpt}),
\[
	\E[\xi\log\xi]<m_1\quad\text{if and only if }\quad\Lambda=\{\g>1:\E|X|^\g<+\infty\text{ and }m_\g<\g m_1\}\neq\emptyset.
\]
Indeed, by Lemma \ref{lemma:derivadamtheta}, $\E[\xi\log\xi]<m_1$ is equivalent to $\Lambda_0=\{\g>1:m_\g<\g m_1\}\neq\emptyset$. Furthermore, by convexity of $r\mapsto m_r/r$ (c.f. ibid.), $\g\in\Lambda_0$ is equivalent to $(1,\g]\subset\Lambda_0$. Then, $\Lambda=\Lambda_0\cap\{\g>1:\E|X|^\g<+\infty\}$ by recalling the hypothesis $X,\xi\in L^\alpha(\p)$, for some $\alpha>1$.
\end{remark}


\subsection{A remarkable martingale}\label{subsection:convM}

In this subsection we analyze a martingale that arises by adapting the ones studied by Bercu in \cite{Bercu} and Bertenghi and Rosales-Ortiz in \cite{MarcoAle} to the setting of the RWES. Its convergence plays the key role in determining the convergence of $\tS$ when $p\in(0,1)$. Let $M=\{M_n\}_{n\geq 1}$ be the process defined as 
\[
	M_n=\frac{(n-1)!}{\Gamma(n+pm_1)}(\tS_{n}-\E\tS_{n}),\quad n\geq1. 
\]

\begin{lemma}\label{lemma:Mesmtg} Fix $p\in(0,1)$. The process $M$ is a martingale with respect to the natural filtration $\{\F_n\}_{n\geq 1}$ of $(X,\varepsilon,U,\xi)$. If $X,\xi\in L^\alpha(\p)$ for some $\alpha>1$, then $M\subset L^\alpha(\p)$.
\end{lemma}
\vspace{-.5cm}\begin{proof}[\bf\color{blue4}\textit{Proof.}] It suffices to show the martingale property. By construction of the increments $\tX$, 
\begin{equation} \label{eq:eqauxcero}
	\E[\tX_{n+1}|\F_n]=(1-p)\E X+\frac{pm_1}{n}\tS_{n},\quad\text{for every }n\geq 1.
\end{equation}
By (\ref{eq:recursionmomentos}), we have that
\[
	\frac{n}{n+p\mxi}\E\left[\tS_{n+1}-\E\tS_{n+1}\Big|\F_n\right]=\tS_{n}-\E\tS_{n}.
\]
Multiplying both sides of the last display by $(n-1)!/\Gamma(n+p\mxi)$, we conclude $\E[ M_{n+1}|\F_n]= M_n$. Secondly, by means of Hölder's inequality, $|\tS_{n}|^\alpha\leq n^{\alpha-1}\sum_{k=1}^{n}|\tX_k|^\alpha$. Letting $q=\alpha$ in Lemma \ref{lemma:recursionmedias}, we conclude that $M\subset L^\alpha(\p)$ whenever if $X,\xi\in L^\alpha(\p)$.
\end{proof}

In sharp contrast to the martingale $M$ in \cite{Bercu} and \cite{MarcoAle} (in which particularly $|\xi|\equiv1$), the conditions $X,\xi\in L^2(\p)$ and $p\in(1/2,1]$ are not sufficient to ensure $L^2(\p)$-boundedness of $M$. As will be shown below, the uniform integrability of $M$ is bound to the condition 
\[
	\E\left[\xi^{1/[1\wedge (pm_1)]}\log\xi\right]<\left[1\wedge (pm_1)\right]m_{1/[1\wedge (pm_1)]},
\]
in deep connection to the condition of convergence to non-degenerate limits for the weights of the memory sub-trees $||T_r(\bullet)||$. 

We remark that the proofs for the strong convergences asserted in this work differ from the supercritical case $pm_1>1$ to the sub-critical and critical cases $pm_1\leq1$. The first one, being the simplest, is basically derived from the results found in Section \ref{section:casep=1} for the pure echoing case $p=1$ and its link with the martingale $M$ of Lemma \ref{lemma:Mesmtg} through the memory sub-trees representation (\ref{eq:Snarbolitos}) found in Section \ref{section:trees}. On the other hand, in the subcritical and critical regimes $pm_1\leq 1$, more care is required both showing the uniform integrability of $M$ and identifying the limiting random variable. Namely, even though we have characterized the convergence of each weighted memory sub-tree $||T_r(\bullet)||$ of the RWES in Corollary \ref{cor:convTreesRootedatk}, further inspection of the martingale $M$ is required to guarantee its overall convergence. Nevertheless, the theory already built in Section \ref{section:trees} is the path to identify the random limit of $M$ in the regime $pm_1\leq1$ when the spins $X$ are centered. Before jumping onto the proof of the main results, which are tackled in the next section, we provide two properties of $M$ that apply to the three regimes, but shall find its use particularly for the case $pm_1\leq1$.


\begin{lemma}\label{lemma:incrsMgamma} Write $\Delta M_n=M_n-M_{n-1}$ for $n\geq2$. If $pm_\alpha\neq1$, there exist constants $C_\alpha,C_\alpha'>0$ such that
\[
	\E|\Delta M_n|^\alpha\leq C_\alpha n^{-\alpha[1\wedge (pm_1)]}\log(n)^{\alpha\1_{pm_1=1}}+C_\alpha'n^{[1\vee (pm_\alpha)]-1-\alpha pm_1}\quad\text{for every }n\geq2.
\]
\end{lemma}
\vspace{-1cm}\begin{proof}[\bf\color{blue4}\textit{Proof.}]  By definition of $M$,
\[
	\frac{\Gamma(n+1+pm_1)}{n!}\Delta M_{n+1}=-\left[\left(\frac{pm_1}{n}\tS_{n}+\E\tX_{n+1}\right)-\left(\tX_{n+1}+\frac{pm_1}{n}\E\tS_{n}\right)\right].
\]
Assume first that $\alpha>2$. Letting $\delta>0$, $w=\frac{pm_1}{n}\tS_{n}+\E\tX_{n+1}$ and $z=\tX_{n+1}+\frac{pm_1}{n}\E\tS_{n}$ in Lemma \ref{lemma:inequalities} gives
\[\begin{split}
	\frac{\Gamma(n+1+pm_1)^\alpha}{n!^\alpha}|\Delta M_{n+1}|^\alpha&\leq (1+\delta N_{\alpha,\delta})\Big|\frac{pm_1}{n}\tS_{n}+\E\tX_{n+1}\Big|^\alpha+N_{\alpha,\delta}'\Big|\tX_{n+1}+\frac{pm_1}{n}\E\tS_{n}\Big|^\alpha+\\
	&-\alpha\Big|\frac{pm_1}{n}\tS_{n}+\E\tX_{n+1}\Big|^{\alpha-2}\left(\frac{pm_1}{n}\tS_{n}+\E\tX_{n+1}\right)\left(\tX_{n+1}+\frac{pm_1}{n}\E\tS_{n}\right),
\end{split}\]
where $N_{\alpha,\delta}\downarrow0$ as $\delta\to0^+$. If $\alpha>2$, choose $\delta=\delta(\alpha)>0$ sufficiently small such that $\delta N_{\alpha,\delta}<\alpha-1$, whereas if $\alpha\in(1,2]$ the expression in the last display still holds removing the constant $\delta N_{\alpha,\delta}$ (c.f. Lemma \ref{lemma:inequalities}). Taking conditional expectation with respect to $\F_{n}$ and noting that
\[
	\E\left[\tX_{n+1}+\frac{pm_1}{n}\E\tS_{n}\Big|\F_{n}\right]=\frac{pm_1}{n}\tS_{n}+\E\tX_{n+1},
\]
by the choice of $\delta>0$, we have
\begin{equation}\label{eq:ineqDeltaMgamma}\begin{split}
	\frac{\Gamma(n+1+pm_1)^\alpha}{n!^\alpha}\E[|\Delta M_{n+1}|^\alpha|\F_{n}]&\leq N_{\alpha,\delta}'\E\left[\Big|\tX_{n+1}+\frac{pm_1}{n}\E\tS_{n}\Big|^\alpha\Big|\F_{n}\right]\\
	&\leq N_{\alpha,\delta}'2^{\alpha-1}\left[(1-p)\E |X|^\alpha+\frac{pm_\alpha}{n}\sum_{i=1}^{n}|\tX_i|^\alpha+\Big|\frac{pm_1}{n}\E\tS_{n}\Big|^\alpha\right],
\end{split}\end{equation}
where (\ref{eq:eqauxcero}) and the elementary inequality $|a+b|^\alpha\leq2^{\alpha-1}(|a|^\alpha+|b|^\alpha)$ were used for the second inequality. The assertion follows by taking expectation, by the asymptotic (\ref{eq:stirling}) and letting $q=\alpha$ in Lemma \ref{lemma:recursionmedias}.
\end{proof}


\vspace{-5mm}The key for proving the convergences in the subcritical and critical regimes $pm_1\leq1$ relies on the following estimate for $\E|M|^\alpha$; although we remark that it applies to the supercritical case $pm_1>1$ as well. Write $\tS^{(\alpha)}$ for the RWES with spins distributed as $|X|^\alpha$ and echoing variables distributed as $\xi^\alpha$, and $M^{(\alpha)}$ for the martingale defined in Lemma \ref{lemma:Mesmtg} associated to $\tS^{(\alpha)}$.

\begin{lemma}\label{lemma:extensiongamma>2} Fix $r>0$ and $\alpha>1$ such that $X,\xi\in L^{\alpha r}(\p)$ and $pm_{\alpha r}\neq1$. There exist constants $C(\alpha,r)>0$ and $\varphi(\alpha,r)\geq0$ such that
\[
	\E\left[\sup_{1\leq k\leq n}\big|M_k^{(r)}\big|^{\alpha}\right]\leq C(\alpha,r)n^{\beta^+(\alpha,r)}\log(n)^{\varphi(\alpha,r)},
\]
for every $n\geq2$, where 
\[
	\beta(\alpha,r)=\begin{dcases}
		[1\vee (pm_{\alpha r})]-\alpha pm_r&\text{ if }\alpha\in(1,2],\\
		\alpha\Big[\Big(\mfrac{1}{2}\vee\sup_{\theta\in[2,\alpha]}\mfrac{pm_{r\theta}}{\theta}\Big)-pm_r\Big]&\text{ if }\alpha>2.
\end{dcases}\]
Moreover, if $\frac{1}{2}\vee\frac{1}{\alpha}\vee \frac{pm_{\alpha r}}{\alpha}<pm_r$, then $\varphi(\alpha,r)=0$.
\end{lemma}

The proof of Lemma \ref{lemma:extensiongamma>2} is quite technical and postponed to the Appendix \ref{section:appendixextensionsgamma>2}. We remit ourselves for the moment to sketch the proof, considering $\alpha\in(2^{\ell-1},2^\ell]$. First, if $\ell=1$, by the Burkholder-Davis-Gundy inequality, 
\begin{equation}\label{eq:sketchextgamma>2}
	\E\left[\sup_{1\leq k\leq n}\big|M_k^{(r)}\big|^{\alpha}\right]\leq C\E\left[[M^{(r)}]_n^{\alpha/2}\right]\leq C\sum_{k=2}^n\E\big|\Delta M_k^{(r)}\big|^{\alpha},
\end{equation}
for some $C>0$ (since $\alpha/2\in(1/2,1]$) and the claim follows by Lemma \ref{lemma:incrsMgamma}. When $\ell=2$, the first inequality in (\ref{eq:sketchextgamma>2}) still holds, but the second one is not immediate, since $\alpha/2>1$ in this case. We thus consider the predictable compensator $A^{(1)}$ of the quadratic variation of $Y^{(0)}:=M^{(r)}$. Applying the Burkholder-Davis-Gundy inequality to the martingale $Y^{(1)}:=[Y^{(0)}]-A^{(1)}$ we find that
\[
	\E\left[[M^{(r)}]_n^{\alpha/2}\right]\leq C'\left(\E\left[[Y^{(1)}]_n^{\alpha/4}\right]+\E[(A_n^{(1)})^{\alpha/2}]\right),
\]
and the result follows by finding an upper bound for $\E[(A_n^{(1)})^{\alpha/2}]$ via the result for $\ell=1$ and by Lemma \ref{lemma:incrsMgamma}, since it is possible to show that
\[
	\E\left[[Y^{(1)}]_n^{\alpha/4}\right]\leq C''\sum_{k=1}^n\E\big|\Delta M_n^{(r)}\big|^{\alpha}.
\]
Inductively, if $\alpha\in(2^{\ell},2^{\ell+1}]$, we consider the predictable compensators $A^{(i)}$ of $[Y^{(i-1)}]$ and apply iteratively the Burkholder-Davis-Gundy inequality to the martingales $Y^{(i)}:=[Y^{(i-1)}]-A^{(i)}$, where $i\in\{1,\dots,\ell\}$. This way, it is possible to write
\[
	\E\left[\sup_{1\leq j\leq n}\big|M_n^{(r)}\big|^{\alpha}\right]\leq C(\alpha,r)\left\{\E\left[[Y^{(\ell)}]_n^{\alpha/2^{\ell+1}}\right]+\sum_{i=1}^\ell\E\left[(A_n^{(i)})^{\alpha/2^i}\right]\right\},
\]
where the estimate for the first addend of the right-hand side is found with Lemma \ref{lemma:incrsMgamma} and the estimates for the second addends are determined by the result being true for powers in $(1,2^{\ell}]$.

\begin{remark}\normalfont The explicit definition of $\varphi(\alpha,r)$ can also be found in Appendix \ref{section:appendixextensionsgamma>2}, equation (\ref{eq:phiexplogaritmos}). We remark, though, that the expression is intricate and the only actual relevant information to our purposes is that it is zero whenever $\frac{1}{2}\vee\frac{1}{\alpha}\vee \frac{pm_{\alpha r}}{\alpha}<pm_r$. Nevertheless, in order to prove Lemma \ref{lemma:extensiongamma>2}, $\varphi(\alpha,r)$ must be specified.
\end{remark}


\subsection{Proof of Theorem \ref{thm:convergencesupercritical} in the supercritical regime $pm_1>1$}\label{subsection:supercriticalcase}

In this subsection we provide proof of Theorem \ref{thm:convergencesupercritical} for the supercritical case $pm_1>1$, which shows that $\tS$ exhibits a super-linear scaling exponent in this regime when $\E[\xi\log\xi]<m_1$, in sharp contrast to the precursors of this model, the step-reinforced and counterbalaced-step random walks, in which the rates of convergence are always less or equal than 1.


\begin{proof}[\bf\color{blue4}\textit{Proof of Theorem ~\ref{thm:convergencesupercritical}}] The case $p=1$ was already proven in Proposition \ref{prop:convergencecasepeq1}, we thus fix $p\in(0,1)$. Recall from (\ref{eq:comportamientoasintoticoesperanzas}) that
\begin{equation}\label{eq:auxsupercritical1}
	n^{-pm_1}\E\tS_n\sim \frac{\E X}{\Gamma(1+pm_1)}\left[1+\frac{1-p}{pm_1-1}\right].
\end{equation}
Since $\E|\tS_n-\E\tS_n|\leq2\E\sum_{r=1}^n|X_r|$ $||T_r(n)||$, substituting $\E|X|$ for $\E X$ in the last display, we find that $M$ is bounded in $L^1(\p)$ and converges almost surely to an integrable limit by Lemma \ref{lemma:Mesmtg} and (\ref{eq:stirling}). Moreover, (\ref{eq:auxsupercritical1}) ensures that $n^{-pm_1}\tS_n$ converges a.s. to an integrable limit. Let us show that such limit is $\sum_{r\geq1}X_r\LL_r^{(p,\xi)}$, where the variables $\{\LL_r^{(p,\xi)}\}_{r\geq1}$ were defined in Corollary \ref{cor:convTreesRootedatk} and are independent of $\{X_r\}_{r\geq1}$.

First of all, the series is well defined, in fact it is absolutely convergent almost surely. Indeed, in virtue of ibid., $\LL_r^{(p,\xi)}=0$ a.s. for every $r\geq1$, if $\E[\xi\log\xi]\geq m_1$. Whereas, if $\E[\xi\log\xi]<m_1$, the variables are non-negative non-degenerate and, by the distributional equation (\ref{eq:eqdistLr}), there exists $C>0$ such that $\E\LL_r^{(p,\xi)}\leq Cr^{-pm_1}$ for every $r\geq1$. Since $pm_1>1$, we have that $\sum_{r\geq1}\E|X_r\LL_r^{(p,\xi)}|<+\infty$ and thus the series is absolutely convergent almost surely.

On the other hand, by means of the representation (\ref{eq:Snarbolitos}) and Corollary \ref{cor:convTreesRootedatk},
\[
	\lim_{n\to\infty}\sum_{r=1}^mX_r\frac{||T_r(n)||}{n^{pm_1}}=\sum_{r=1}^mX_r\LL_r^{(p,\xi)}\quad\text{a.s. for every }m\geq1.
\]
Since we have found that $n^{-pm_1}\tS_n$ converges almost surely,
\[
	\lim_{n\to\infty}\sum_{r=m}^{n}X_r\frac{||T_r(n)||}{n^{pm_1}}=\lim_{n\to\infty}n^{-pm_1}\tS_n-\sum_{r=1}^mX_r\LL_r^{(p,\xi)}\quad\text{a.s.}\quad\text{for every }m\geq1.
\]
Thus, the limit as $m\to\infty$ of the left-hand side of the last display exists a.s. and the assertion follows by proving that 
\[
	\lim_{m\to\infty}\lim_{n\to\infty}\sum_{r=m}^n|X_r|\frac{||T_r(n)||}{n^{pm_1}}=0\quad\text{a.s.}
\]
This is true by an application of Fatou's lemma, since (\ref{eq:esperanzaTrn}) gives
\[
	\sum_{r=m}^n\E\left[|X_r|\text{ }\frac{||T_r(n)||}{n^{pm_1}}\right]\leq C\sum_{r=m}^nr^{-pm_1},\quad\text{for every }n\geq m\geq1.
\]
Therefore, $n^{-pm_1}\tS_n\to\sum_{r\geq1}X_r\LL_r^{(p,\xi)}$ almost surely, where $\sum_{r\geq1}X_r\LL_r^{(p,\xi)}$ is non-degenerate if $\E[\xi\log\xi]< m_1$ and is zero a.s., otherwise.

We are thus only left to show that the convergence holds in $L^\gamma(\p)$ for every $\gamma\in\Lambda$ if $\E[\xi\log\xi]<m_1$. By Remark \ref{remark:Lambdanonempty}, $\Lambda\neq\emptyset$ and $\g\in\Lambda$ iff $(1,\g]\subset\Lambda$. Moreover, since $pm_1>1$, for every $\g\in\Lambda$ we have $\frac{1}{2}\vee\frac{1}{\g}\vee\frac{pm_\g}{\g}<pm_1$. Thus, by letting $r=1$ and $\alpha=\g$ in Lemma \ref{lemma:extensiongamma>2}, $\varphi(\g,1)=\beta^+(\g,1)=0$ and there exists $C(\g)>0$ (that does not depend on $n$) such that $\E|M_n|^\gamma\leq C(\g)$ for every $n\geq1$. Thus, $M$ is bounded in $L^\gamma(\p)$ and converges in such sense.
\end{proof}

Remark that one can also show that $M$ is bounded in $L^\g(\p)$ in the supercritical case when $\E[\xi\log\xi]<m_1$, without Lemma \ref{lemma:extensiongamma>2}. Namely, it is possible to construct an estimate for $\E|M_n|^\g$ by making use of $n^{-pm_1}||T_r(n)||\to\LL_r^{(p,\xi)}$ a.s. and in $L^\g(\p)$ (c.f. Corollary \ref{cor:convTreesRootedatk}ii)) and writing $\E[(\LL_r^{(p,\xi)})^\g]$ in virtue of the distributional equation (\ref{eq:eqdistLr}). We close this subsection by pointing out that, with the proof of Theorem \ref{thm:convergencesupercritical}, it was also shown that

\begin{corollary}\label{cor:limitseriessupercritical} Let $p\in(0,1)$. If $pm_1>1$, then 
\[
	\lim_{n\to\infty}\frac{\tS_{n}}{n^{pm_1}}=\sum_{r\geq1}X_r\LL_r^{(p,\xi)}\quad\text{a.s.},
\]
where $\{\LL_r^{(p,\xi)}\}_{r\geq1}\subset L^1(\p)$ are the random variables defined in Corollary \ref{cor:convTreesRootedatk}, they are independent of $\{X_r\}_{r\geq1}$ and the series is absolutely convergent a.s. and integrable. Furthermore, if $\E[\xi\log\xi]< m_1$, the series is non-degenerate, whereas otherwise, $\sum_{r\geq1}X_r\LL_r^{(p,\xi)}=0$ almost surely.
\end{corollary}


\subsection{Proof of Theorem \ref{thm:convergencesubcritical}}\label{subsection:subcriticalcase}

The goal of this subsection is to provide proof of Theorem \ref{thm:convergencesubcritical} for the subcritical $pm_1<1$ and critical $pm_1=1$ regimes. Remark that the convergences presented in this section broadly generalize that of the (positively) step-reinforced random walks ($\xi\equiv1$), specifically in Theorem 1 in \cite{NRBM}. Indeed, in such setting $\xi\equiv1$, so $m_1=1$ and $pm_1<1$ since $p$ is a probability parameter. Since $m_2=1$, the condition $m_2<2m_1$ is naturally met. So the convergences in $L^2(\p)$ and almost surely stated in Theorem 1 in ibid., respectively, follow by Theorem \ref{thm:convergencesubcritical}.

First, we prove the strong convergence of the rescaling of the RWES around its mean for the critical and sub-critical cases $pm_1\leq1$ towards a non-degenerate random variable. Subsequently, a corollary regarding the distribution of the limit is provided by means of the link of the RWES with URRT in a simile fashion as that of the supercritical regime $pm_1>1$. Lastly, we end the subsection by characterizing the hypothesis in Theorem \ref{thm:convergencesubcritical}ii).


\begin{proof}[\bf\color{blue4}\textit{Proof of Theorem~\ref{thm:convergencesubcritical}}] By Remark \ref{remark:Lambdanonempty}, $\Lambda\neq\emptyset$. Then, it suffices to show that $M$ is bounded in $L^{\g}(\p)$ for $\g\in\Lambda$. In such case, in conjunction with Lemma \ref{lemma:recursionmedias} and (\ref{eq:stirling}), we can infer that the limits
\[\begin{dcases}
	\lim_{n\to\infty}n^{-pm_1}\left[\tS_n-\mfrac{(1-p)\E X}{1-pm_1}\text{ }n\right]&\text{if }pm_1\in(1/2,1),\\
	\lim_{n\to\infty}n^{-1}\left[\tS_n-(1-p)\E X\text{ }n\log(n)\right]&\text{if }pm_1=1,
\end{dcases}\]
exist both almost surely and in $L^\g(\p)$, and are thus non-degenerate and centered. To this end, we claim that the constants in Lemma \ref{lemma:extensiongamma>2} are $\varphi(\g,1)=\beta^+(\g,1)=0$ for every $\gamma\in\Lambda\cap(1/pm_1,+\infty)$. Indeed, if $pm_1=1$ then $\Lambda\cap(1/pm_1,+\infty)=\Lambda$ and the condition
\[
	\frac{1}{2}\vee\frac{1}{\g}\vee \frac{pm_\g}{\g}< pm_1
\]
is met trivially. For the subcritical regime $pm_1<1$, this condition is met for every $\g\in\Lambda\cap(1/pm_1,+\infty)$, where such intersection is non-empty by hypothesis. By letting $r=1$ and $\alpha=\g\in\Lambda\cap(1/pm_1,+\infty)$ in Lemma \ref{lemma:extensiongamma>2}, we thus have that $\varphi(\g,1)=0$ and there exists $C(\g)>0$ such that
\[
	\E\left[\sup_{1\leq k\leq n}\big|M_k\big|^{\g}\right]\leq C(\g)n^{\beta^+(\g,1)} \quad\text{for every }n\geq2,
\]
where
\[
	\beta(\g,1)=\begin{dcases}
		[1\vee (pm_\g)]-\g pm_1&\text{ if }\g\in(1,2],\\
		\g\Big[\Big(\frac{1}{2}\vee\sup_{\theta\in[2,\g]}\left\{\frac{pm_\theta}{\theta}\Big\}\Big)-pm_1\right]&\text{ if }\g>2,
	\end{dcases}\leq0,
\]
since $\g\in\Lambda$ iff $(1,\g]\subset\Lambda$ (c.f. Remark \ref{remark:Lambdanonempty}). Hence, $M$ is bounded in $L^\gamma(\p)$ and the proof is complete.
\end{proof}


In virtue of the representation (\ref{eq:Snarbolitos}) of $\tS$ in terms of its memory sub-trees, it is possible to retrieve a representation of the limiting random variable when the spins of the RWES are centered, $\E X=0$. Moreover, in general, we can also derive information about the distribution of the limit in terms of the partial sums of the weighted memory sub-trees as follows.

\begin{corollary}\label{cor:seriespm1geq1} Let $pm_1\in(1/2,1]$ and consider $\{\LL_{r}^{(p,\xi)}\}_{r\geq1}$, the random variables defined in Corollary \ref{cor:convTreesRootedatk}. Under the conditions of Theorem \ref{thm:convergencesubcritical}, the following assertions hold:
\vspace{-0.5cm}\begin{enumerate}[leftmargin=0.5cm]
	\item[\textbf{{\color{blue4}{i)}}}] If $\E X=0$, then
	\[
		\lim_{n\to\infty}n^{-pm_1}\tS_n=\sum_{r\geq1}X_r\LL_r^{(p,\xi)}\quad\text{a.s.},
	\]
	where the series is non-degenerate, absolutely convergent a.s., and $\{\LL_r^{(p,\xi)}\}_{r\geq1}\indep\{X_r\}_{r\geq1}$.
	\item[\textbf{{\color{blue4}{ii)}}}] There exists $\gamma\in(1,2]$ such that
	\[\begin{dcases}
		\lim_{n\to\infty}n^{-pm_1}\left[\tS_n-\mfrac{(1-p_)\E X}{1-pm_1}\text{ }n\right]=\lim_{m\to\infty}\sum_{r=1}^m\left(X_r\LL_r^{(p,\xi)}-(1-p\1_{r>1})\mfrac{\E X(r-1)!}{\Gamma(r+pm_1)}\right)&\text{if }pm_1<1,\\
		\lim_{n\to\infty}n^{-pm_1}\left[\tS_n-(1-p)\E X\text{ }n\log(n)\right]=\lim_{m\to\infty}\sum_{r=1}^m\left(X_r\LL_r^{(p,\xi)}-(1-p\1_{r>1})\mfrac{\E X}{r}\right)&\text{if }pm_1=1,
	\end{dcases}\]
	hold in $L^\gamma(\p)$. In particular, the limits in the left-hand side are $F_\infty$-distributed, where $F_\infty$ is the limiting distribution of the partial sums in the right-hand side as $m$ tends to $\infty$.
\end{enumerate}
\end{corollary}
\begin{proof}[\bf\color{blue4}\textit{Proof.}] Write $M_\infty=\lim_{n\to\infty}n^{-pm_1}\tS_n$. By Theorem \ref{thm:convergencesubcritical}, the limit of the martingale $M$ of Lemma \ref{lemma:Mesmtg} and $M_\infty$ exist a.s. and coincide pointwise. In addition, the limits hold in $L^\gamma(\p)$ for every $\g\in\Lambda$ as well. Recall also that $\LL_r^{(p,\xi)}$ is the limit of $\TT_{r,n}:=\mfrac{(n-1)!}{\Gamma(n+pm_1)}\E||T_r(n)||$, both a.s. and in $L^\gamma(\p)$ for every $r\geq1$, where $\{T_r(n)\}_{n\geq r\geq1}$ is independent of $\{X_r\}_{r\geq1}$ in virtue of representation (\ref{eq:Snarbolitos}).

For i), similarly to the proof of Theorem 1 in \cite{NRBM}, if the spins $X$ are centered, the process $\sum_{r=1}^mX_r\LL_r^{(p,\xi)}$ is a martingale with respect to the filtration $\{\mathcal{G}_m\}_m\geq1$, where $\mathcal{G}_m=\s(\{\LL_r^{(p,\xi)}\}_{r\geq1},\{X_r\}_1^m)$. Without loss of generality, assume that $\gamma\in(1,2]$. By letting $w=\sum_{r=1}^{k-1}X_r\LL_r^{(p,\xi)}$ and $z=X_k\LL_k^{(p,\xi)}$ in Lemma \ref{lemma:inequalities} and taking conditional expectation with respect to $\mathcal{G}_{k-1}$ iteratively for $k\in\{2,\dots,m\}$, one gets
\[
	\E\Bigg|\sum_{r=1}^mX_r\LL_r^{(p,\xi)}\Bigg|^\gamma\leq N\E|X|^\gamma\sum_{r=1}^m\E[(\LL_r^{(p,\xi)})^\gamma],
\]
where $N>0$ does not depend on any parameter. We claim that therefore $\sum_{r=1}^mX_r\LL_r^{(p,\xi)}$ is bounded in $L^\gamma(\p)$ and the series $\sum_{r\geq1}X_r\LL_r^{(p,\xi)}$ exists a.s. as the limit of this martingale. Indeed, since $\LL_1^{(p,\xi)}\in L^\gamma(\p)$, by the distributional identity (\ref{eq:eqdistLr}) there exists a constant $C_\gamma>0$ such that $\E[(\LL_r^{(p,\xi)})^\gamma]\leq C_\gamma r^{-\gamma pm_1}$ for every $r\geq1$. Hence, by ii) and completeness of $L^\gamma(\p)$, we infer that $M_\infty=\sum_{r\geq1}X_r\LL_r^{(p,\xi)}$ almost surely.

In order to prove ii) note that by (\ref{eq:Snarbolitos}), the right-hand side of ii) equals
\[
	\sum_{r=1}^m\left(X_r\LL_r^{(p,\xi)}-\E X\E\LL_r^{(p,\xi)}\right)=\lim_{n\to\infty}\sum_{r=1}^m\left(X_r\TT_{r,n}-\E X\E\TT_{r,n}\right) \quad\text{a.s.}.
\]
Recall from (\ref{eq:Tr(n)esmtg}) that $\{\TT_{r,n}\}_{n\geq r}$ is a uniformly integrable martingale with respect to $\{\F_n\}_{n\geq r}$. Thus, $\E\TT_{r,n}=\E\LL_r^{(p,\xi)}$ and, by representation (\ref{eq:Snarbolitos}), for every $m\geq1$,
\[
	\lim_{n\to\infty}\sum_{r=m+1}^n\left(X_r\TT_{r,n}-\E X\E\TT_{r,n}\right)=M_\infty-\sum_{r=1}^m\left(X_r\LL_r^{(p,\xi)}-\E X\E\LL_r^{(p,\xi)}\right) \quad\text{a.s.}
\]
In addition, by construction,
\[
	\sum_{r=m+1}^n\left(X_r\TT_{r,n}-\E X\E\TT_{r,n}\right)=(M_n-M_m)-\sum_{i=m+1}^n\sum_{r=1}^mX_r\Delta\TT_{r,i}\quad\text{for every }n\geq m\geq1, 
\]
so, for every $m\geq1$,
\[
	\sum_{r=1}^m\left(X_r\LL_r^{(p,\xi)}-\E X\E\LL_r^{(p,\xi)}\right)=M_m+\lim_{n\to\infty}\sum_{i=m+1}^n\sum_{r=1}^mX_r\Delta\TT_{r,i}\quad\text{a.s.}
\]
Since $M_m\to M_\infty$ in $L^\gamma(\p)$, ii) therefore folllows by showing that
\begin{equation}\label{eq:aux2seriespmleq1}
	\lim_{m\to\infty}\E\Bigg|\lim_{n\to\infty}\sum_{i=m+1}^n\sum_{r=1}^mX_r\Delta\TT_{r,i}\Bigg|^\gamma=0.
\end{equation}
In the setting of Lemma \ref{lemma:inequalities}, by iteratively letting $w=\sum_{i=m+1}^{k}\sum_{r=1}^mX_r\Delta\TT_{r,i}$ and $z=\sum_{r=1}^mX_r\Delta\TT_{r,k+1}$, and taking conditional expectation with respect to $\F_{k}$ for $k\in\{m+1,\dots,n-1\}$, one gets
\begin{equation}\label{eq:aux3seriespmleq1}
	\E\Bigg|\sum_{i=m+1}^n\sum_{r=1}^mX_r\Delta\TT_{r,i}\Bigg|^\gamma\leq N\sum_{i=m+1}^n\E\Bigg|\sum_{r=1}^mX_r\Delta\TT_{r,i}\Bigg|^\gamma.
\end{equation}
To estimate the expectations in the right-hand side of the last display, note that by construction
\[
	\sum_{r=1}^mX_r\Delta\TT_{r,i}=-\frac{pm_1}{i-1+pm_1}\sum_{r=1}^mX_r\TT_{r,i-1}+\frac{(i-1)!}{\Gamma(i+pm_1)}\varepsilon_i\xi_i\tX_{\U[i]}\sum_{r=1}^m\1_{\U[i]\in T_r(i-1)}.
\]
Then, by letting $w$ and $z$ be the first and second addends of the last display in Lemma \ref{lemma:inequalities}, respectively, and taking conditional expectation with respect to $\F_{i-1}$,
\begin{equation}\label{eq:aux4seriespmleq1}
	\E\Bigg|\sum_{r=1}^mX_r\Delta\TT_{r,i}\Bigg|^\gamma\leq N\E|X|^\g\frac{pm_\gamma}{i-1}\frac{(i-1)!^\g}{\Gamma(i+pm_1)^\gamma}\sum_{r=1}^m\E\left[\sum_{j\in T_r(i-1)}\omega(j)^\gamma\right],
\end{equation}
where the definition of $\omega(\cdot)$ can be recalled from (\ref{eq:defwtree}). By plugging (\ref{eq:aux4seriespmleq1}) in (\ref{eq:aux3seriespmleq1}) and applying (\ref{eq:esperanzaTrn}) and (\ref{eq:stirling}),
\[
	\E\Bigg|\sum_{i=m+1}^n\sum_{r=1}^mX_r\Delta\TT_{r,i}\Bigg|^\gamma\leq C_\g\left(\sum_{r=1}^m\frac{(r-1)!}{\Gamma(r+pm_\gamma)}\right)\left(\sum_{i=m+1}^n\frac{\Gamma(i+pm_\gamma)}{\Gamma(i+1+\gamma pm_1)}\right).
\]
Lastly, by the gamma identities (\ref{eq:sumagammas}) and (\ref{eq:stirling}) and Fatou's Lemma,
\[
	\E\Bigg|\lim_{n\to\infty}\sum_{i=m+1}^n\sum_{r=1}^mX_r\Delta\TT_{r,i}\Bigg|^\gamma\leq C_\g'\lim_{n\to\infty} (m^{1\vee pm_\g-\g pm_1}+n^{pm_\gamma-\gamma pm_1}),
\]
from where (\ref{eq:aux2seriespmleq1}) follows.
\end{proof}


To conclude, we give a characterization of the condition for convergence in the subcritical case $pm_1<1$ stated in Theorem \ref{thm:convergencesubcritical}ii), so as a sufficient and a necessary condition for it to hold.

\begin{proposition} Let $pm_1\in(0,1)$. The following assertions hold:
\vspace{-0.5cm}\begin{enumerate}[leftmargin=0.7cm]
	\item[\textbf{{\color{blue4}{i)}}}] $\Lambda\cap(1/pm_1,+\infty)\neq\emptyset$ is equivalent to $1/pm_1\in\Int\{a>0:X,\xi\in L^{a}(\p)\}$ and $pm_{1/pm_1}<1$.
	\vspace{-2mm}\item[\textbf{{\color{blue4}{ii)}}}] If $X,\xi\in L^{1/pm_1}(\p)$ and $\E[\xi^{1/pm_1}\log\xi]< m_{1/pm_1}pm_1$, then i) holds.
	\vspace{-2mm}\item[\textbf{{\color{blue4}{iii)}}}] If $\Lambda\cap(1/pm_1,+\infty)\neq\emptyset$, then $pm_1\in(1/2,1)$.
\end{enumerate}
\end{proposition}\vspace{-.5cm}
\begin{proof}[\bf\color{blue4}\textit{Proof.}] Let $\phi_1(r)=m_r/r$ for $r\in\I$, where $\I=\Int\{a>0:x,\xi\in L^a(\p)\}$. By Lemma \ref{lemma:derivadamtheta}, $\phi_1$ is continuously differentiable and strictly convex on $\I$. To prove i), assume first that $\Lambda\cap(1/pm_1,+\infty)\neq\emptyset$. Since there exists $\g>1/pm_1$ such that $m_\g<\g m_1$, then $1/pm_1\in\I$. By convexity, $\phi_1(1/pm_1)<m_1$, which in turn is equivalent to $pm_{1/pm_1}<1$. Now, supposing that $1/pm_1\in\I$ and $\phi_1(1/pm_1)<m_1$, by continuity there exists $\delta>0$ such that $\g:=1/pm_1+\delta\in\I$ and $\phi_1(\g)<m_1$. Thus, $1\vee pm_{\g}<\g pm_1$ and i) is proven. 

On the other hand, assume that $X,\xi\in L^{1/pm_1}(\p)$. If the condition in ii) is met, by Lemma \ref{lemma:derivadamtheta}, $\phi_1$ is strictly decreasing in $[1,1/pm_1+\delta]$ for some $\delta>0$. Thus, there exists some $\gamma>1/pm_1$ such that $\phi_1(\g)<m_1$. Hence, $1\vee pm_\g<\g pm_1$ and $\Lambda\cap(1/pm_1,+\infty)\neq\emptyset$. 

Lastly, we address that iii) follows by construction.
\end{proof}


\subsection{Proof of Laws of Large Numbers for the subcritical and critical regimes $pm_1\leq1$}\label{subsection:LLN}

We conclude the current section by giving proof of the Law of Large Numbers for the RWES asserted in Theorem \ref{thm:LLN}. Remark that this result generalizes Theorem 1.1 in \cite{MarcoAle}, given for the case $\xi\equiv1$. In this direction, a modification of the second Law of Large Numbers for martingales will be required. Substituting $M_n^2$ with $|M_n|^\gamma$ in the proof of Theorem 1.3.17 in \cite{RandIterModels}, one readily gets, 


\begin{proposition}\label{prop:2LLNMtgs} Let $\{M_n\}_{n\geq1}$ be a martingale such that $M\subset L^\alpha(\p)$ for some $\alpha>1$. If there exists $\beta_0>0$ such that the sequence $n^{-\beta_0}\E|M_n|^\alpha$ is bounded, then for every $\beta>\beta_0/\alpha$,
\[
	\lim_{n\to\infty}\frac{M_n}{n^{\beta}}=0\quad\text{almost surely and in }L^\alpha(\p).
\]
\end{proposition}
\vspace{-0.5cm}Sketching the proof in ibid., in virtue of Doob's maximal inequality, it is plain to see that for every $n\geq0$ and $\delta>0$,
\[
	\p\left(\sup_{2^{n}\leq k\leq2^{n+1}}\frac{|M_n|}{n^{\beta}}>\delta\right)\leq\frac{\E|M_{2^{n+1}}|^\alpha}{\delta^\alpha 2^{n\alpha\beta}}\leq2^{n(\beta_0-\alpha\beta)}\left(\frac{\alpha}{\alpha-1}\right)^\alpha\frac{2^{\beta_0}}{\delta^\alpha }\sup_{k\geq1}\left\{\frac{\E|M_{k}|^\alpha}{k^{\beta_0}}\right\},
\]
where the supremum in the right-hand side of the last display is finite by hypothesis. Thus, if $\beta>\beta_0/\alpha$,
\[
	\sum_{n\geq0}\p\left(\sup_{2^{n}\leq k\leq2^{n+1}}\frac{|M_n|}{n^{\beta}}>\delta\right)<+\infty,
\]
and the almost sure converge holds in virtue of Borel-Cantelli lemma. To argue for the $L^\alpha(\p)$ convergence, we simply remark that, 
\[
	\lim_{n\to\infty}\frac{\E|M_n|^\alpha}{n^{\alpha\beta}}\leq\lim_{n\to\infty} n^{\beta_0-\alpha\beta}\sup_{k\geq1}\left\{\frac{\E|M_{k}|^\alpha}{k^{\beta_0}}\right\}=0.
\]


\begin{proof}[\bf\color{blue4}\textit{Proof of Theorem~\ref{thm:LLN}}] The result i) for $pm_1=1$ readily follows from Lemma \ref{lemma:recursionmedias} and Theorem \ref{thm:convergencesubcritical}. To prove ii), fix $pm_1\in(0,1)$. We aim to show that
\begin{equation}\label{eq:aux1LLN}
	\text{there exist}\quad \beta_0>0,\text{ }\gamma>1\quad\text{such that}\quad\frac{\beta_0}{\g}<1-pm_1\quad\text{and}\quad \E|M_n|^\g=\text{o}(n^{\beta_0}).
\end{equation}
This way, in virtue of Proposition \ref{prop:2LLNMtgs}, Lemma \ref{lemma:recursionmedias} and (\ref{eq:stirling}),
\[
	\lim_{n\to\infty}\frac{1}{n}\left[\tS_n-\frac{(1-p)\E X}{1-pm_1}\text{ }n\right]=\lim_{n\to\infty}\frac{M_n}{n^{1-pm_1}}=0\quad\text{a.s. and in }L^\g(\p),
\]
which gives the desired convergence. By Lemma \ref{lemma:extensiongamma>2}, for every $\alpha>1$ such that $X,\xi\in L^\alpha(\p)$ and $pm_\alpha\neq1$, we have
\[
	\E\left[\sup_{1\leq k\leq n}\big|M_k\big|^\alpha\right]\leq C(\alpha) n^{\beta^+(\alpha,1)}\log(n)^{\varphi(\alpha,1)},
\]
for some $C(\alpha)>0$ and $\beta^+(\alpha,1),\varphi(\alpha,1)\geq0$. In order to prove (\ref{eq:aux1LLN}), we shall argue that there exists $\gamma>1$ such that the lemma can be applied, such that
\[
	\beta^+(\g,1)=\begin{dcases}
		\left[\mfrac{\g}{2}\vee1\vee(pm_\g)\right]-\g pm_1&\text{ if }\E[\xi\log\xi]\geq m_1,\\
		\left[\left(\mfrac{\g}{2}\vee1\right)-\g pm_1\right]^+&\text{ if }\E[\xi\log\xi]< m_1,
	\end{dcases}
\]
where $\g$ may be different from  line to line, and such that (\ref{eq:aux1LLN}) holds by letting
\begin{equation}\label{eq:aux2LLN}
	\beta_0=\begin{cases}
		\frac{1}{2}\left[\left(\frac{\g}{2}\vee1\vee (pm_\g)\right)+\g\right]-\g pm_1&\text{ if }\E[\xi\log\xi]\geq m_1,\\
		\frac{1}{2}\left[\left(\frac{\g}{2}\vee1\right)+\g\right]-\g pm_1&\text{ if }\E[\xi\log\xi]< m_1\text{ and }pm_1<\frac{1}{2}\vee\frac{1}{\g},\\
		1-pm_1&\text{ otherwise}.
	\end{cases}
\end{equation}

In the case $\E[\xi\log\xi]\geq m_1$, by Lemma \ref{lemma:derivadamtheta}, $r\mapsto m_r/r$ is strictly increasing in $\Int\{a\geq1:X,\xi\in L^a(\p)\}$. Then, by continuity, $\gamma>1$ can be chosen such that $pm_\g\neq1$ and $pm_\g/\g<1$. Hence, $\beta^+(\g,1)=\frac{\g}{2}\vee1\vee(pm_\g)-\g pm_1$. By construction, $\frac{\g}{2}\vee1\vee (pm_\g)<\g$, so by writing $\beta_0$ as in (\ref{eq:aux2LLN}), we see that $\beta^+(\g,1)<\beta_0<\g(1-pm_1)$. Thus, (\ref{eq:aux1LLN}) holds and so the stated convergence.

Suppose now that $\E[\xi\log\xi]<m_1$. By Remark \ref{remark:Lambdanonempty}, $\Lambda\neq\emptyset$ and $\g\in\Lambda$ iff $(1,\g]\subset\Lambda$. Then, $\beta^+(\g,1)=(\frac{\g}{2}\vee1-\g pm_1)^+$ for every $\g\in\Lambda$. If $\beta^+(\g,1)>0$, since $\frac{\g}{2}\vee1<\g$, by defining $\beta_0$ as in (\ref{eq:aux2LLN}) gives $\beta^+(\g,1)<\beta_0<\g(1-pm_1)$ and (\ref{eq:aux1LLN}) holds. Lastly, if $\beta^+(\g,1)=0$, then $\E|M_n|^\g=\OO(\log(n)^{\varphi(\g,1)})$. Thus (\ref{eq:aux1LLN}) is true in particular by writing $\beta_0=1-pm_1$.
\end{proof}


\section{Distributional properties of the limiting random variables}\label{section:propsLo}

The last part of this work is devoted to give some distributional properties of the limiting random variables $\{\LL_r^{(p,\xi)}\}_{r\geq1}$ involved in the convergence of the weighted memory sub-trees $||T_r(n)||$ of $\tS_n=\sum_{r=1}^nX_r||T_r(n)||$ (see Corollary \ref{cor:convTreesRootedatk}). In virtue of the distributional equation of Lemma \ref{lemma:fixedpointeq||T||} for the case of pure echoing $p=1$, it is possible to derive one associated to limiting random variable $\LL^{(\xi)}$ defined in Proposition \ref{prop:convergencecasepeq1}. Jointly with our findings in ibid., this provides the means to discover the following properties of the distribution of $\LL^{(\xi)}$, arguing similarly as in \cite{Lucille} for the ERW. It should be highlighted that the proof of the assertions (ii)-(vi) involving the distributional equation and its consequences are greatly inspired in \cite{Lucille}, which is based on similar methods exposed in \cite{Chauvin}.


\begin{proposition}\label{prop:propsLo} Suppose that $\E[\xi\log\xi]< m_1$. The random variable $\LL^{(\xi)}$ defined in Proposition \ref{prop:convergencecasepeq1} has the following properties:\vspace{-.5cm}
\begin{enumerate}[leftmargin=0.7cm]
	\item[\textbf{{\color{blue4}{i)}}}] Let $p_0=\p(\xi=0)$ and $\xi^+$ have distribution $\p(\xi^+\in\dd x)=\p(\xi\in\dd x|\xi>0)$. If $p_0>0$, then 
	\[
		\LL^{(\xi)}\myeqd\LL^{(\xi^+)}\M_{1-p_0}^{m_1/(1-p_0)},
	\]
	where $\M_{1-p_0}$ is a Type II Mittag-Leffler$(1-p_0)$ random variable, independent of $\LL^{(\xi^+)}$.
	\item[\textbf{{\color{blue4}{ii)}}}] $\LL^{(\xi)}$ satisfies the distributional equation
	\begin{equation}\label{eq:fixedpointeqLo}
		\LL^{(\xi)}\myeqd V^{m_1}\hat{\LL}^{(\xi)}+\xi(1-V)^{m_1}\check{\LL}^{(\xi)},
	\end{equation}
	where $\hat{\LL}^{(\xi)},\check{\LL}^{(\xi)}$ are independent copies of $\LL^{(\xi)}$, $V$ is a standard uniform random variable and all these random elements are mutually independent.
	\item[\textbf{{\color{blue4}{iii)}}}] The mean of $\LL^{(\xi)}$ is $\Gamma(1+m_1)^{-1}$. The condition $m_k<km_1$ for $k\geq2$ is equivalent to $\LL^{(\xi)}\in L^k(\p)$. In such case, we have the recursive relation
	\begin{equation}\label{eq:mukLo}
		\E[(\LL^{(\xi)})^k]=\frac{\Gamma(km_1+1)^{-1}}{km_1-m_k}\sum_{i=1}^{k-1}{k\choose i }m_i\Gamma(1+im_1)\Gamma(1+(k-i)m_1)\E[(\LL^{(\xi)})^i]\E[(\LL^{(\xi)})^{k-i}].
	\end{equation}
	\item[\textbf{{\color{blue4}{iv)}}}] The variance $\Var(\LL^{(\xi)})=0$ if and only if $\xi=1$ almost surely. Furthermore, whenever $\p(\xi=1)<1$ the random variable $\LL^{(\xi)}$ is supported on $[0,\infty)$, while $\LL^{(\xi)}=1$ if $\xi=1$ almost surely.
	\item[\textbf{{\color{blue4}{v)}}}] The characteristic function $\varphi$ of $\LL^{(\xi)}$ satisfies the relation
	\begin{equation}\label{eq:igualdadvarphiVns}
		\varphi(t)=\E\left[\varphi(tV^{m_1})\varphi(t\xi[1-V]^{m_1})\right]\quad\text{for every }t\in\R,
	\end{equation}
	where $V$ is a standard uniform random variable independent of $\xi$. Furthermore, $\LL^{(\xi)}$ is non-atomic if $\p(\xi=1)<1$.
	
	\item[\textbf{{\color{blue4}{vi)}}}] Assume that $\p(\xi=1)<1$ and that $\xi>0$ a.s. If there exists a constant $\alpha\in(0,1/m_1)$ such that $\xi^{-\alpha}\in L^1(\p)$, then $|\varphi(t)|=\OO(|t|^{-\alpha})$. If, in addition, $\E\xi<1$ and such $\alpha>1$, then $\LL^{(\xi)}$ has a bounded continuous density. On the other hand, whenever $\xi^{-1}\in L^\infty(\p)$ for every $k\geq1$,
	\[
		|\varphi(t)|=\mathcal{O}(|t|^{-k})\quad\text{as }|t|\to\infty
	\]
	 and, thus, $\LL^{(\xi)}$ has a smooth, bounded and $\Leb$-almost surely positive density and admits moments of all orders.
\end{enumerate}
\end{proposition}


Due to the detailed description of the distribution of the limiting random variable $\LL^{(\xi)}$ for the case of pure echoing $p=1$ provided in Proposition \ref{prop:propsLo}, one can derive the following properties for the distribution of $\LL_r^{(p,\xi)}$ for $p\in(0,1)$ directly from the equality in distribution (\ref{eq:eqdistLr}). We only address that in virtue of Proposition \ref{prop:propsLo}i), we have 
\[
	\LL_1^{(p,\xi)}\myeqd\LL^{(\xi^+)}\M_{p(1-p_0)}^{m_1/(1-p_0)},
\]
where $\M_{p(1-p_0)}\sim$Mittag-Leffler$(p(1-p_0))$ is independent of $\LL^{(\xi^+)}$ and $p_0=\p(\xi=0)$. This way, the existence of a density for $\LL_1^{(p,\xi)}$ stated in Proposition \ref{prop:distribucionLr(p)}v) relies only on the condition $1/\xi^+\in L^\infty(\p)$, by means of Proposition \ref{prop:propsLo}vi).

\begin{proposition}\label{prop:distribucionLr(p)} Let $p\in(0,1]$ and suppose that such that $\E[\xi\log\xi]<m_1$. The random variables $\{\LL_r^{(p,\xi)}\}_{r\geq1}$ defined in Corollary \ref{cor:convTreesRootedatk} exhibit the following properties:\vspace{-0.5cm}
\begin{enumerate}[leftmargin=0.7cm]
	\item[\textbf{{\color{blue4}{i)}}}] For every $r\geq1$, the variable $\LL_r^{(p,\xi)}$ is supported on $[0,\infty)$ and $\Var\LL_r^{(p,\xi)}>0$.
	\item[\textbf{{\color{blue4}{ii)}}}] The mean of $\LL_r^{(p,\xi)}$ is $\Gamma(1+pm_1)^{-1}$ if $r=1$, and $(1-p)\frac{(r-1)!}{\Gamma(r+pm_1)}$ otherwise. 
	\item[\textbf{{\color{blue4}{iii)}}}] Let $k\geq2$. The condition $m_k<km_1$ is equivalent to $\{\LL_r^{(p,\xi)}\}_{r\geq1}\subset L^k(\p)$. In such case,
\[
	\E[(\LL_r^{(p,\xi)})^i]=(1-p\1_{r>1})(r-1)!\frac{\Gamma(1+ipm_1)}{\Gamma(r+ipm_1)}\E[(\LL^{(\varepsilon\xi)})^i]\quad(1\leq i\leq k).
\]
	\item[\textbf{{\color{blue4}{iv)}}}] Let $\xi^+$ have the distribution $\p(\xi^+\in\dd x)=\p(\xi\in\dd x|\xi>0)$. Whenever $1/\xi^+\in L^\infty(\p)$:
	\begin{itemize}[leftmargin=0.5cm, label=\textcolor{blue4}{$\blacktriangleright$}]
		\item $\LL_1^{(p,\xi)}$ is absolutely continuous with respect to the Lebesgue measure,
		\item If $r\geq2$, then $\LL_r^{(p,\xi)}$ has an only atom at $0$ with mass $p$ and its Lebesgue decomposition is thus given by $\delta_{\{0\}}$ and $(1-p)\p(\LL_1^{(p,\xi)}\beta_r^{pm_1}\in\dd x)$.
	\end{itemize} 
\end{enumerate}
\end{proposition}

By last, we provide proof for Proposition \ref{prop:propsLo} and highlight once more that the properties (ii)-(vi) are obtained by adapting the arguments found in \cite{Lucille} to the RWES. We shall therefore remit to ibid. for details when possible.


 \begin{proof}[\bf\color{blue4}\textit{Proof of Proposition~\ref{prop:propsLo}}] 
 
\noindent {\color{blue4}{\textit{\textbf{i)}}}} In order to show the assertion i), it is necessary to go back to the context of the continuous time BRW associated to $\tS$ when $p=1$ (c.f. Subsection \ref{subsection:BRW}) . If $p_0=\p(\xi=0)>0$, denoting $\xi^+$ as the law of $\xi$ conditioned on $\xi>0$ and $\delta_0\sim$Bernoulli$(1-p_0)$, we have $\xi\myeqd\delta_0\xi^+$. Write $\ZZ^{(\xi)}$ and $\ZZ^{(\xi^+)}$ for the BRWs associated to echo variables $\xi$ and $\xi^+$, respectively. 

Let $(Y_t^{(p)})_{t\geq0}$ denote a Yule process with rate $p$, i.e., a pure-birth process starting at 1 with birth rate from state $n\in\N$ given by $pn$. It is well known that $(Y_t^{(p)})_{t\geq0}\myeqd(Y_{pt}^{(1)})_{t\geq0}$. Then, by construction $\ZZ_t^{(\xi)}\myeqd\ZZ_{(1-p_0)t}^{(\xi^+)}$, which in turn implies that $\W_\infty^{(\xi)}\myeqd\W_\infty^{(\xi^+)}$. From (\ref{eq:limSncomoWyE}) we remark that $\LL^{(\xi)}$ is a function of the echo variables $\{\xi_n\}_{n\geq1}$ and the genealogical tree $\T$. Then, it is independent of the Yule process $(Y_t)_{t\geq0}$ and, thus, $\LL^{(\xi)}\indep\mathcal{E}_\infty$. This argument also applies to $\LL^{(\xi^+)}$ and the limit $\mathcal{E}_\infty^{(\xi^+)}$ related to its Yule process $Y^{(\xi^+)}$. Then,
\[
	\E\e^{i\theta\log\LL^{(\xi)}}\Gamma(1+i\theta m_1)=\E\e^{i\theta\log\W_\infty^{(\xi)}}=\E\e^{i\theta\log\W_\infty^{(\xi^+)}}=\E\e^{i\theta\log\LL^{(\xi^+)}}\Gamma\left(1+i\theta \mfrac{m_1}{1-p_0}\right).
\]
for every $\theta\in\R$. The assertion readily follows by justifying that if $\M_{1-p_0}$ is a Type II Mittag-Leffler random variable with parameter $1-p_0$, then
\[
	\E\exp\left\{i\theta\log(\M_{1-p_0}^{m_1/(1-p_0)})\right\}=\frac{\Gamma\left(1+i\theta \mfrac{m_1}{1-p_0}\right)}{\Gamma(1+i\theta m_1)}.
\]

In general, if $\M\sim$Mittag-Leffler($q$) with $q\in(0,1)$, then $\M$ is positive almost surely and its Mellin transform is given by $\E\M^a=\Gamma(1+a)/\Gamma(1+qa)$ for $a>-1$ \cite{ML}. Thus, the moment generating function of $\log\M$ is finite in the non-degenerate interval $(-1,+\infty)$. Then, $z\mapsto\E\M^z$ is analytic in the strip $D=\{z\in\C:\text{Re}z>-1\}$ (c.f. Exercise 4.8 \cite{Kallenberg}) and coincides in $(-1,+\infty)$ with the mapping $z\mapsto\Gamma(1+z)\Gamma(1+qz)^{-1}$, also analytic in $D$. Thus, $z\mapsto\E\M^z$ equals $\Gamma(1+z)\Gamma(1+qz)^{-1}$ in $D$ by the Identity Theorem. 
 
 
\noindent {\color{blue4}{\textit{\textbf{ii)}}}} According to Lemma \ref{lemma:fixedpointeq||T||}, for every $n\geq1$ there is the equality in distribution $\tS_{n+1}\myeqd \hatS_{R_n}+\xi\checkS_{n+1-R_n}$; where $ \hatS,\checkS$ are copies of $\tS$, $R_n$ denotes the number of red balls after $n-1$ extractions in a Pólya-Eggenberger urn with initial condition $(1,1)$, and all these elements are mutually independent. Therefore
\[
	\lim_{n\to\infty}\frac{\tS_{n+1}}{(n+1)^{m_1}}\myeqd\lim_{n\to\infty}\left[\left(\frac{R_n}{n+1}\right)^{m_1}\frac{\hatS_{R_n}}{R_n^{m_1}}+\xi\left(1-\frac{R_n}{n+1}\right)^{m_1}\frac{\checkS_{n+1-R_n}}{(n+1-R_n)^{m_1}}\right],
\]
which gives the distributional equation in ii) since $n^{-1}R_n\to V$ a.s. as $n\to\infty$, where $V$ is a standard uniform random variable (c.f. Theorem 3.2 in \cite{Mahmoud}).


\noindent {\color{blue4}{\textit{\textbf{iii)}}}} From Proposition \ref{prop:convergencecasepeq1} we already know that $\E\LL^{(\xi)}=\Gamma(1+m_1)^{-1}$ if $\E[\xi\log\xi]<m_1$ and that $\LL^{(\xi)}\in L^k(\p)$ if $m_k<km_1$. To show that $\LL^{(\xi)}\in L^k(\p)$ implies that $m_k<km_1$, by means of (\ref{eq:fixedpointeqLo}),
\begin{equation}\label{eq:mukLoparcial}
	\E[(\LL^{(\xi)})^k]=\frac{1+m_k}{1+km_1}\E[(\LL^{(\xi)})^k]+\sum_{i=1}^{k-1}{k\choose i}m_i\beta(1+im_1,1+(k-i)m_1)\E[(\LL^{(\xi)})^i]\E[(\LL^{(\xi)})^{k-i}].
\end{equation} 
From the last display it is easy to see that, assuming $m_k>km_1$ leads to $\E[(\LL^{(\xi)})^k]<0$, and assuming $m_k=km_1$ leads to $\LL^{(\xi)}=0$ a.s., both contradictions. Therefore it must be that $m_k<km_1$ and (\ref{eq:mukLo}) follows from (\ref{eq:mukLoparcial}).


\noindent {\color{blue4}{\textit{\textbf{iv)}}}} If $\xi\equiv1$, then by construction $\tS(n)\equiv n+1$ and thus $\LL^{(\xi)}\equiv1$. Secondly, if $\Var\LL^{(\xi)}=0$, then $\LL^{(\xi)}=\Gamma(1+m_1)^{-1}$ a.s. by ii). In virtue of the distributional equation (\ref{eq:fixedpointeqLo}), then $\xi=\frac{1-V^{m_1}}{(1-V)^{m_1}}$ a.s. But $\xi$ and $V$ are independent, so we conclude $\xi\equiv1$. Therefore, $\Var\LL^{(\xi)}=0$ is equivalent to $\xi=1$ a.s. To complete the argument, it remains to show that $\supp\LL^{(\xi)}=[0,\infty)$ if $\p(\xi=1)<1$. But this can be done exactly as in Theorem 2.5 in \cite{Lucille} via the distributional equation (\ref{eq:fixedpointeqLo}). Namely, taking elements in $\supp\LL^{(\xi)}$ around 1, instead of around 0, and separating into the cases $m_1<1$, $m_1=1$ and $m_1>1$.


\noindent {\color{blue4}{\textit{\textbf{v)}}}} Conditioning on $(V,\xi)$ we get (\ref{eq:igualdadvarphiVns}) from (\ref{eq:fixedpointeqLo}). Following \cite{Lucille}, in addition this gives
\begin{equation}\label{eq:ineqphiLo}
	|\varphi(t)|\leq\E[|\varphi(tV^{m_1})\varphi(t\xi[1-V]^{m_1})|].
\end{equation}
Note that with the elementary bound $|\varphi|\leq1$ we have $|\varphi(t)|\leq\E|\varphi(tV^{m_1})|$ from the latter. And thus, inductively, if $\{V_n\}_{n\geq1}$ are independent standard uniform random variables, conditioning on $V_n,\dots V_2$ successively gives
\begin{equation}\label{eq:desigualdadvarphiVns}
	|\varphi(t)|\leq\E|\varphi(tV_1^{m_1})|\leq\E|\varphi(tV_1^{m_1}V_2^{m_1})|\leq\cdots\leq\E|\varphi(tV_1^{m_1}\cdots V_n^{m_1})|.
\end{equation}
Thus, if $\p(\xi\neq1)>0$ then $\supp\LL^{(\xi)}=[0,\infty)$ and we can conclude that $\LL^{(\xi)}$ is non-atomic by showing that $|\varphi(t)|\to0$ as $|t|\to\infty$. This follows from (\ref{eq:ineqphiLo}) and (\ref{eq:desigualdadvarphiVns}) with the same arguments as in the three-step proof for Theorem 2.5 in \cite{Lucille}. Indeed, since we assumed $\xi>0$ a.s. and there is independence between $\xi(1-V)^{m_1}$ and the independent sequence of standard uniform random variables $\{V_n\}_{n\geq1}$, it is possible to substitute $(1-V)^{m_1}$ with $\xi(1-V)^{m_1}$ in the proof given in ibid. to get the same conclusion. Therefore we refer to \cite{Lucille} for the details.

\noindent {\color{blue4}{\textit{\textbf{vi)}}}} Since $|\varphi(t)|\to0$ as $|t|\to\infty$, for every $\varepsilon>0$ there is a $T_\varepsilon>0$ such that $|\varphi(t)|<\varepsilon$ if $|t|\geq T_\varepsilon$. We proceed to modify slightly the arguments in Theorem 2.5 \cite{Lucille} to adapt it to RWES. Note that
\[
	|\varphi(t)|\leq\varepsilon\E|\varphi(tV^{m_1})|+\p\left(\xi V^{m_1}\leq \frac{T_\varepsilon}{|t|}\right),
\]
in virtue of (\ref{eq:ineqphiLo}). Thus, provided $\xi^{-\alpha}\in L^1(\p)$ for some $\alpha\in(0,m_1^{-1})$, the last expression gives
\[
	|\varphi(t)|\leq\varepsilon\E|\varphi(tV^{m_1})|+\left(\frac{T_\varepsilon}{|t|}\right)^\alpha\frac{m_{-\alpha}}{1-\alpha m_1},\quad|t|\geq T_\varepsilon,
\]
applying Markov's inequality to the integrable random variable $(\xi V^{m_1})^{-\alpha}$. Thus, by applying this inequality inductively with (\ref{eq:ineqphiLo}) and (\ref{eq:desigualdadvarphiVns}), one gets
\begin{equation}\label{eq:ineqaux2cor4652}
	|\varphi(t)|\leq\left(\frac{T_\varepsilon}{|t|}\right)^\alpha\frac{m_{-\alpha}}{1-\alpha m_1-\varepsilon},
\end{equation}
for every $\varepsilon\in(0,1-\alpha m_1)$ and $|t|\geq T_\varepsilon$. Therefore, $|\varphi(t)|=\OO(|t|^{-\alpha})$. If we in addition consider the hypotheses $m_1<1$ and $\alpha>1$, then $\LL^{(\xi)}$ has a continuous and bounded Leb-almost surely density, by Fourier's inversion Theorem.

On the other hand, by assuming that $\p(\xi\geq z_0)=1$ for some $z_0>0$, from (\ref{eq:ineqaux2cor4652}) follows that
\begin{equation}\label{eq:ineqauxcor465}
	|\varphi(t)|\leq\varepsilon\E|\varphi(tV^{m_1})|+\frac{m_{-\alpha}}{1-\alpha m_1-\varepsilon}\left(\frac{T_\varepsilon}{|t|}\right)^\alpha\E\left[ V^{-1}\1_{V\geq1-\left(\frac{T_\varepsilon}{|t|z_0}\right)^{\alpha}}\right].
\end{equation}
Then following the proof of Theorem 2.5 in \cite{Lucille}, this gives $|\varphi(t)|=\OO(|t|^{-2\alpha})$ and, repeating the argument inductively yields $|\varphi(t)|=\OO(|t|^{-k})$, for every $k\geq1$. Therefore, $\LL^{(\xi)}$ has a smooth and bounded density and the proof of (vi) is complete.
\end{proof}


\section*{Acknowledgements} I would like to thank Jean Bertoin and Alejandro Rosales-Ortiz for our always stimulating discussions, their comprehensive review on this work and extensive suggestions to improve it. In particular, to Jean Bertoin for introducing me to the world of random walks with step-reinforcement and sharing his brilliant intuition. Also, to Alejandro Rosales-Ortiz for his careful reading of the proofs and suggestions, particularly for the section devoted to Branching Random Walks. I would also like to thank Justinas Grigaitis for our discussions about noisy memory models in behavioral economics.


\appendix

\section{Proof of the lemmas for the pure echoing case}\label{section:appendixUchiyama}


We provide proof for Lemmas \ref{lemma:derivadamtheta} and \ref{lemma:inequalities} in subsection \ref{subsection:||T||acotLa}. It is worth pointing out that the main argument was taken from \cite{Uchiyama} and that the proof is included here for sake of completeness (since our hypotheses are slightly different).

\begin{proof}[\bf\color{blue4}\textit{Proof of Lemma~\ref{lemma:derivadamtheta}}] Recall the notation $m_\alpha=\E\xi^\alpha$. In virtue of $\xi\in L^\alpha(\p)$ for $\alpha>1$, it is possible to show with standard dominated convergence arguments that $\theta\mapsto m_\theta$ is twice differentiable on $\I$, strictly convex and $\frac{\dd}{\dd \theta} m_\theta=\E[\xi^\theta\log\xi]$. Then, for every fixed $\theta\in\I$ the mapping $\phi_\theta:\R^+\to\R^+$ defined as $\phi_\theta(r)= m_{r\theta}/r$ is continuous and twice differentiable in $\I$. To see that $\phi_\theta$ is strictly convex in $\I$, compute
\begin{equation}\label{eq:derivadaphi}
	\phi'_\theta(r)=\frac{1}{r}\left[\frac{\dd}{\dd r}m_{\theta r}-\phi_\theta(r)\right]\quad\text{and}\quad\phi''_\theta(r)=\frac{1}{r}\left[\frac{\dd^2}{\dd r^2}m_{\theta r}-2\phi'_\theta(r)\right].
\end{equation}
Since $\phi_\theta(r)\to+\infty$ as $r\to0^+$, there is a positive neighborhood of 0 such that $\phi'_\theta$ is strictly negative. If this holds for all $r\in\I$, $\phi_\theta$ is obviously convex and strictly decreasing. If there is a point in $\I$ in which $\phi'_\theta$ is non-negative, by the intermediate value theorem there exists $r_\theta\in\I$ such that $\phi'_\theta(r_\theta)=0$. The last display yields then that $\phi''_\theta(r_\theta)>0$, since $\alpha\mapsto m_\alpha$ is strictly convex. Furthermore $\phi_\theta$ attains its minimum in $r_\theta$ and only in $r_\theta$ due to the strict convexity of $\alpha\mapsto m_\alpha$.

From the former analysis we conclude that $\phi$ is strictly convex in $\I$ and its limit as $r\to\infty$ exists, although it may be infinite. In any case, for every $\theta\in \I$ there is $r_\theta\in(0,+\infty]$ such that
\[
	\lim_{r\to r_\theta^-}\phi_\theta(r)=\inf_{r\in\R^+}\phi_\theta(r).
\]
Since $\phi_\theta$ is strictly decreasing in $(0,r_\theta)$, the condition $\E[\xi^\theta\log\xi]<\E\xi^\theta/\theta$ is equivalent to $\phi_\theta'(1)<0$. This again happens if and only if $1<r_\theta$ and then for every $r\in(1,r_\theta)$ we have $m_{\theta r}/r\theta<m_\theta/\theta$. Lastly, remark that $r_\theta<+\infty$ is equivalent to $\phi_\theta(r_\theta)=\frac{\dd}{\dd r}m_{\theta r}|_{r=r_\theta}$ by (\ref{eq:derivadaphi}), completing the proof.
\end{proof}

Moving over to the proof of Lemma \ref{lemma:inequalities}, we point out that $\gamma\in(1,2]$ is a particular case of inequality (4.2) in \cite{Uchiyama}.

\begin{proof}[\bf\color{blue4}\textit{Proof of Lemma~\ref{lemma:inequalities}}] For $\gamma>2$, similarly to ibid., writing $|1+x|^\gamma=[(1+x)^2]^{\gamma/2}$ and approximating by Taylor series, we have that $|1+x|^\gamma-1-\gamma x=$o$(x)$ as $x\to0$. Also, $|1+x|^\gamma-1-\gamma x\sim |x|^{\gamma}$ as $x\to\pm\infty$. Thus, for $\delta>0$ we define the finite quantities
\[
	N_{\gamma,\delta}=\sup_{x'\in[-\delta,\delta]}\frac{|1+x'|^\gamma-1-\gamma x'}{|x'|}\quad\text{and}\quad N_{\gamma,\delta}'=\sup_{|x'|>\delta}\frac{|1+x'|^\gamma-1-\gamma x'}{|x'|^\gamma}
\]
and note that for every $x\in\R$,
\begin{equation}\label{eq:auxmidesigUchiyama}
	|1+x|^\gamma\leq\begin{cases}
		1+\gamma x+N_{\gamma,\delta}|x|&\text{ if }|x|\leq\delta,\\
		1+\gamma x+N_{\gamma,\delta}'|x|^\gamma&\text{ if }|x|>\delta.
	\end{cases}
\end{equation}
Furthermore, by a standard analysis of its derivative, we find that $x\mapsto\frac{|1+x|^\gamma-1-\gamma x}{|x|}$ is decreasing in $[-\delta,0)$ and increasing in $(0,\delta]$. Thus, $N_{\gamma,\delta}\downarrow0$ as $\delta\to0^+$. In a similar manner, $x\mapsto\frac{|1+x|^\gamma-1-\gamma x}{|x|^\gamma}$ has an asymptote in 0, is increasing in $(-\infty,-\delta]$ and decreasing in $[\delta,\infty)$. Then, $N_{\gamma,\delta}'\downarrow1$ as $\delta\to\infty$. Lastly, the inequalities in the lemma follow by letting $x=z/w$ in (\ref{eq:auxmidesigUchiyama}) for $w\neq0$; and trivially for the case $w=0$, since $N_{\gamma,\delta}'>1$.
\end{proof}


\section{Proof of Lemma \ref{lemma:extensiongamma>2}}\label{section:appendixextensionsgamma>2}

This last section is devoted to prove Lemma \ref{lemma:extensiongamma>2}, which gives an estimate of $\E[\sup_{1\leq k\leq n}|M_n^{(r)}|^{\alpha}]$. In that direction, let us introduce some notation. For $\alpha>1$, $r>0$, define
\[
	\varphi_1(\alpha,r)=\1\left\{pm_r=\mfrac{1}{\alpha}\right\}+\1\left\{pm_r=\mfrac{1}{\alpha}\vee\mfrac{pm_{\alpha r}}{\alpha}\right\}.
\]
In addition, if $\alpha\in(2^{\ell-1},2^\ell]$ for some $\ell\geq2$, define recursively
\[\begin{split}
	\varphi_m(\alpha,r)&=\varphi_1(\alpha,r)+\sum_{i=1}^{m-1}\1\Big\{pm_r\leq\mfrac{1}{2^{i+1}}\vee\mfrac{1}{\alpha}\vee\sup_{\theta\in[2^i,\alpha]}\mfrac{pm_{r\theta}}{\theta}\Big\}\text{ }\varphi_{m-i}\left(\mfrac{\alpha}{2^i},2^ir\right)+\\
		&+\sum_{i=1}^{m-1}\frac{\alpha}{2^i}\Big[\varphi_1(2^i,r)+\1\left\{pm_r\leq \mfrac{1}{2^i}\vee\mfrac{ pm_{2^ir}}{2^i}\right\}\1_{pm_{2^ir}=1}+\1\Big\{pm_r=\mfrac{1}{2^{i+1}}\vee\mfrac{1}{\alpha}\vee\sup_{\theta\in[2^i,\alpha]}\mfrac{pm_{r\theta}}{\theta}\Big\}\Big],
\end{split}\]
for $m\in\{2,\dots,\ell\}$. Then, we write
\begin{equation}\label{eq:phiexplogaritmos}
	\varphi(\alpha,r)=\varphi_{\lceil\log\alpha/\log2\rceil}\left(\alpha,r\right)\quad\text{and}\quad\beta(\alpha,r)=\begin{dcases}
		[1\vee (pm_{\alpha r})]-\alpha pm_r&\text{ if }\alpha\in(1,2],\\
		\alpha\Big[\Big(\mfrac{1}{2}\vee\sup_{\theta\in[2,\alpha]}\mfrac{pm_{r\theta}}{\theta}\Big) -pm_r\Big]&\text{ if }\alpha>2.
	\end{dcases}
\end{equation}

Recall that $\tS^{(r)}$ denotes the RWES with spins distributed as $|X|^r$ and echoing variables distributed as $\xi^r$. Denote by $M^{(r)}$ the martingale defined in Lemma \ref{lemma:Mesmtg} associated to $\tS^{(r)}$, provided that $X,\xi\in L^r(\p)$. Write also $[Y]_n=\sum_{k=2}^n(\Delta Y_n)^2$ for the quadratic variation of a process $\{Y_n\}_{n\geq1}$.


\begin{proof}[\bf\color{blue4}\textit{Proof of Lemma~\ref{lemma:extensiongamma>2}}] Let $\ell:=\lceil\log\alpha/\log2\rceil$ and note that $\alpha\in(2^{\ell-1},2^\ell]$. Remark that if $\frac{1}{2}\vee\frac{1}{\alpha}\vee\frac{pm_{\alpha r}}{\alpha}<pm_r$, Lemma \ref{lemma:derivadamtheta} guarantees that $pm_{r\theta}/\theta<pm_r$ for every $\theta\in(1,\alpha]$. Thus, in this case, $\varphi(\alpha,r)=0$ by construction, as stated in the proposition. We now give proof for the inequality by induction over $\ell\geq1$. 

\noindent{\color{blue4}{$\blacktriangleright$}} For the initial case, fix $\ell=1$, so $\alpha\in(1,2]$. Assume that $X,\xi\in L^{\alpha r}(\p)$ and that $pm_{\alpha r}\neq1$. In virtue of the Burkholder-Davis-Gundy inequality, for some $C'(\alpha,r)>0$,
\begin{equation}\label{eq:extensiongammaaux1}
	\E\left[\sup_{1\leq k\leq n}\big|M_k^{(r)}\big|^{\alpha}\right]\leq C'(\alpha,r)\E\left[[M^{(r)}]_n^{\alpha/2}\right]\leq C'(\alpha,r)\sum_{k=2}^n\E\big|\Delta M_k^{(r)}\big|^{\alpha},
\end{equation}
where the second inequality holds since $\alpha/2\in(1/2,1]$. In virtue of Lemma \ref{lemma:incrsMgamma},
\[
	\sum_{k=2}^n\E\big|\Delta M_k^{(r)}\big|^{\alpha}\leq C''(\alpha,r)\left[ \sum_{k=2}^nk^{-\alpha(1\wedge pm_r)}\log(k)^{\alpha\1_{pm_r=1}}+\sum_{k=2}^nk^{[1\vee (pm_{\alpha r})]-1-\alpha pm_r}\right].
\]
Analysing the cases suggested by the minimum and maximum of the last expression, we find that there exist constants $c(\alpha,r),c'(\alpha,r)>0$ such that for every $n\geq1$
\begin{equation}\label{eq:sumaminimos}
	\begin{dcases}
	\sum_{k=2}^nk^{-\alpha(1\wedge pm_r)}\log(k)^{\alpha\1_{pm_r=1}}\leq c(\alpha,r)n^{(1-\alpha pm_r)^+}\log(n)^{\1_{pm_r=\frac{1}{\alpha}}}.\\
	\sum_{k=2}^nk^{[1\vee (pm_{\alpha r})]-1-\alpha pm_r}\leq c'(\alpha,r) n^{ [\{1\vee (pm_{\alpha r})\}-\alpha pm_r]^+}\log(n)^{\1_{pm_r=\frac{1}{\alpha}\vee\frac{pm_{\alpha r}}{\alpha}}}.
	\end{dcases}
\end{equation}
Hence, for some $C(\alpha,r)>0$,
\[\begin{split}
	\E\left[\sup_{1\leq k\leq n}\big|M_k^{(r)}\big|^{\alpha}\right]&\leq C(\alpha,r)n^{\left[\{1\vee (pm_{\alpha r})\}-\alpha pm_r\right]^+}\log(n)^{ \1_{pm_r=\frac{1}{\alpha}}+\1_{pm_r=\frac{1}{\alpha}\vee\frac{pm_{\alpha r}}{\alpha}} }.
\end{split}\]
The result follows by noting that the power of $n$ in the last display is $\beta^+(\alpha,r)$, and that of $\log(n)$ is $\varphi(\alpha,r)=\varphi_1(\alpha,r)$.


\noindent{\color{blue4}{$\blacktriangleright$}} For the inductive step, fix $\ell\geq2$ and suppose that for every $\ell'\in\{1,\cdots,\ell\}$, if $r>0$ and $\alpha\in(2^{\ell'-1},2^{\ell'}]$ are such that $pm_{\alpha r}\neq1$ and $X,\xi\in L^{\alpha r}(\p)$, then there exists $C(\alpha,r)>0$ such that 
\[
	\E\left[\sup_{1\leq k\leq n}\big|M_k^{(r)}\big|^{\alpha}\right]\leq C(\alpha,r)n^{\beta^+(\alpha,r)}\log(n)^{\varphi_{\ell'}(\alpha,r)},
\]
for every $n\geq2$, where $\beta$ and $\varphi_{\ell'}$ were defined in (\ref{eq:phiexplogaritmos}). Let $\gamma\in(2^{\ell},2^{\ell+1}]$ such that $X,\xi\in L^{\g r}(\p)$ and $pm_{\g r}\neq1$. In order to prove the assertion for $\gamma$, for $i\in\{0,\dots,\ell\}$ define $\{A_n^{(i)}\}_{n\geq1}$ and $\{Y_n^{(i)}\}_{r\geq1}$ recursively as
\[\begin{dcases}
	Y^{(0)}=M^{(r)},&\quad A^{(0)}\equiv0,\\
	A_1^{(i)}=0,&\quad A_n^{(i)}=\sum_{k=2}^n\E\left[\big(\Delta Y_k^{(i-1)}\big)^2\big|\F_{k-1}\right],\\
	Y_1^{(i)}=0,&\quad Y_n^{(i)}=[Y^{(i-1)}]_n-A_n^{(i)}=\sum_{k=2}^n\left\{\big(\Delta Y_k^{(i-1)}\big)^2-\E\left[\big(\Delta Y_k^{(i-1)}\big)^2\big|\F_{k-1}\right]\right\}.
\end{dcases}\]
Note that for every $i\in\{1\dots,\ell\}$, $A^{(i)}$ is the predictable compensator of the quadratic variation $[Y^{(i-1)}]$, thus $Y^{(i)}=[Y^{(i-1)}]-A^{(i)}$ is a martingale. From this moment on throughout the proof we will always write $C$ for positive constants following an inequality, these may change from line to line, but they only depend on $\gamma$, $r$ or $i\in\{1,\cdots,\ell\}$. Applying the Burkholder-Davis-Gundy inequality twice,
\[\begin{split}
	\E\left[\sup_{1\leq k\leq n}\big|M_n^{(r)}\big|^{\gamma}\right]&\leq C\E\left[[Y^{(0)}]_n^{\gamma/2}\right]\leq C\left(\E\left[[Y^{(1)}]_n^{\gamma/4}\right]+\E\left[\big(A_n^{(1)}\big)^{\gamma/2}\right]\right).
\end{split}\]
Thus, inductively we get
\begin{equation}\label{eq:auxextgamma1}
	\E\left[\sup_{1\leq k\leq n}\big|M_n^{(r)}\big|^{\gamma}\right]\leq C\left(\E\left[[Y^{(\ell)}]_n^{\gamma/2^{\ell+1}}\right]+\sum_{i=1}^\ell\E\left[\big(A_n^{(i)}\big)^{\gamma/2^i }\right]\right),
\end{equation}
for some positive constant $C>0$.

First, $\gamma/2^{\ell+1}\in(0,1)$ by hypothesis and, thus,
\[
	\E\left[[Y^{(\ell)}]_n^{\gamma/2^{\ell+1}}\right]\leq\sum_{k=2}^n\E\left[\big|\Delta Y_k^{(\ell)}\big|^{\g/2^\ell }\right]\leq\left(\prod_{j=1}^\ell 2^{\gamma/2^j}\right)\sum_{k=2}^n\E\big|\Delta M_k^{(r)}\big|^{\g},
\]
where we have used the definition of $Y^{(i)}$ and Jensen's inequality iteratively for the last inequality. Then, similarly to the initial step $\ell=1$, by Lemma \ref{lemma:incrsMgamma}, there is a constant $C>0$ such that
\begin{equation}\label{eq:auxextgamma2}
	\E\left[[Y^{(\ell)}]_n^{\gamma/2^{\ell+1}}\right]\leq C n^{ \left[\{1\vee (pm_{\gamma r})\}-\g pm_r\right]^+}\log(n)^{\varphi_1(\g,r)}.
\end{equation}

Secondly, we find estimates for $\E[(A_n^{(i)})^{\gamma/2^i}]$. Once again, by applying Jensen's inequality inductively one can write
\[
	\E\left[\big(\Delta Y_k^{(i-1)}\big)^2\big|\F_{k-1}\right]\leq2\times2\E\left[\big(\Delta Y_k^{(i-2)}\big)^4\big|\F_{k-1}\right]\leq\cdots\leq\left(\prod_{j=1}^{i-1}2^{2^{j}}\right)\E\left[\big(\Delta M_k^{(r)}\big)^{2^i}\big|\F_{k-1}\right].
\]
Then, applying (\ref{eq:stirling}) and Lemma \ref{lemma:recursionmedias} in (\ref{eq:ineqDeltaMgamma}), 
\[
	\E\left[\big(\Delta M_k^{(r)}\big)^{2^i}\big|\F_{k-1}\right]\leq C\left( k^{-2^ipm_r}+k^{-2^i[1\wedge (pm_r)]}\log(k)^{2^i\1_{pm_r=1}}+ k^{-2^ipm_r-1}\tS_{k-1}^{(2^ir)}\right).
\]
Thus, by construction of $A^{(i)}$,
\begin{equation}\label{eq:auxextgamma3}\begin{split}
	\big|\big|A_n^{(i)}\big|\big|_{\gamma/2^i}&\leq C \sum_{k=1}^n\Big|\Big|\E\left[\big(\Delta M_k^{(r)}\big)^{2^i}\Big|\F_{k-1}\right]\Big|\Big|_{\gamma/2^i}\\
	&\leq C\left(\sum_{k=2}^n k^{-2^ipm_r}+\sum_{k=2}^nk^{-2^i[1\wedge (pm_r)]}\log(k)^{2^i\1_{pm_r=1}}+\sum_{k=2}^nk^{-2^ipm_r-1}\big|\big|\tS_{k-1}^{(2^ir)}\big|\big|_{\gamma/2^i}\right).
\end{split}\end{equation}
Now we apply the induction hypothesis. Note that for every $i\in\{1,\cdots,\ell\}$, we have $\gamma/2^i\in(2^{\ell-i},2^{\ell-i+1}]$, $pm_{\frac{\gamma}{2^i} 2^ir}=pm_{\gamma r}\neq1$ and $X,\xi\in L^{\frac{\gamma}{2^i}2^ir}(\p)$. Therefore, the bound of the induction hypothesis applies to $M^{(2^ir)}$. Since $\tS_k^{(2^ir)}=\frac{\Gamma(k+pm_{2^ir})}{(k-1)!}M^{(2^ir)}_k+\E\tS_k^{(2^ir)}$, letting $q=2^ir$ in Lemma \ref{lemma:recursionmedias}, 
\[
	\big|\big|\tS_{k-1}^{(2^ir)}\big|\big|_{\gamma/2^i}\leq C\left( k^{pm_{2^ir}+\frac{2^i}{\g}\beta^+\left(\frac{\g}{2^i},2^ir\right)}\log(k)^{\frac{2^i}{\gamma}\varphi_{\ell+1-i}\left(\frac{\g}{2^i},2^ir\right)}+k^{1\vee (pm_{2^ir})}\log(k)^{\1_{pm_{2^ir}=1}}\right).
\]
for every $i\in\{1,\dots,\ell\}$. Thus, by plugging this expression in (\ref{eq:auxextgamma3}), noting that
\[
	\frac{2^i}{\g}\beta\left(\frac{\g}{2^i},2^ir\right)=\begin{dcases}
		\mfrac{2^\ell}{\g}[1\vee(pm_{\g r})]-pm_{2^\ell r}&\text{ if }i=\ell,\\
		\Big(\mfrac{1}{2}\vee\sup_{\theta\in[2^{i+1},\g]}\mfrac{pm_{r\theta}}{\theta/2^i}\Big)-pm_{2^ir}&\text{ if }1\leq i\leq\ell-1,
	\end{dcases}
\]
and by solving the sums with help of (\ref{eq:sumaminimos}), one gets
\[
	\E\left[\big(A_n^{(i)}\big)^{\gamma/2^i}\right]\leq C\left(n^{\big\{\g\big[\big(\frac{1}{2^i}\vee\sup_{\theta\in[2^i,\g]}\frac{pm_{r\theta}}{\theta}\big)-pm_r\big]\big\}^+}\log(n)^{\varphi'_i(\gamma,r)}\right),
\]
where
\[\begin{split}
	\varphi'_i(\gamma,r)&=\1\Big\{pm_r\leq\mfrac{1}{2^{i+1}}\vee\mfrac{1}{\g}\vee\sup_{\theta\in[2^i,\g]}\mfrac{pm_{r\theta}}{\theta}\Big\}\text{ }\varphi_{\ell+1-i}\left(\mfrac{\g}{2^i},2^ir\right)+\\
	&+\mfrac{\g}{2^i}\Big[\varphi_1(2^i,r)+\1\left\{pm_r\leq\mfrac{1}{2^i}\vee\mfrac{ pm_{2^ir}}{2^i}\right\}\1_{pm_{2^ir}=1}+\1 \Big\{pm_r=\mfrac{1}{2^{i+1}}\vee\mfrac{1}{\g}\vee\sup_{\theta\in[2^i,\g]}\mfrac{pm_{r\theta}}{\theta}\Big\}\Big].
\end{split}\]
Recognizing from this expression the addends in the definition of $\varphi_{\ell+1}(\g,r)$, the proof is then completed by substituting the sum over $i\in\{1,\dots,\ell\}$ of the last display added to (\ref{eq:auxextgamma2}) in (\ref{eq:auxextgamma1}) and by noting that
\[
	\max\Big\{[1\vee (pm_{\gamma r})]-\g pm_r,\text{ }\g\Big[\Big(\mfrac{1}{2^i}\vee\sup_{\theta\in[2^i,\g]}\mfrac{pm_{r\theta}}{\theta}\Big)-pm_r\Big]\Big\}\leq\beta(\g,r),
\]
for every $i\in\{1,\dots,\ell\}$.
\end{proof}


\addcontentsline{toc}{section}{References}

\begin{small}

\end{small}


\end{document}